\documentclass[12pt,a4paper]{article}
\usepackage[T2A]{fontenc}
\usepackage[cp1251]{inputenc}
\usepackage[english]{babel}
\usepackage{epsfig,amssymb,epic,eepic,color, graphicx}
\usepackage{fancyhdr}
\usepackage{amsthm}
\usepackage{amsmath}
\usepackage{indentfirst}
\usepackage{setspace}
\usepackage{fancyhdr}
\newtheorem{theorem}{Theorem}

\newtheorem{collary}{Corollary}[section]

\newtheorem{lemma}{Lemma}[section]
\newtheorem{proposition}{Proposition}[section]

\begin{document}
УДК: 517.938.5


{\bf Сharacteristic space of orbits of Morse-Smale diffeomorphisms on surfaces} 

Е.\,V.~Nozdrinova, О.\,V.~Pochinka, Е.\,V. Tsaplina

Nozdrinova Elena Vyacheslavovna, maati@mail.ru Laboratory of Dynamical Systems and Applications, NRU HSE, 25/12, Bolshaya Pecherskaya St., Nizhny Novgorod, Russia, 603155

Pochinka Olga Vitalievna, olga-pochinka@yandex.ru Laboratory of Dynamical Systems and Applications, NRU HSE, 25/12, Bolshaya Pecherskaya St., Nizhny Novgorod, Russia, 603155

Tsaplina Ekaterina Vadimovna, ktsaplina11@mail.ru Laboratory of Dynamical Systems and Applications, NRU HSE, 25/12, Bolshaya Pecherskaya St., Nizhny Novgorod, Russia, 603155

\begin{abstract}
The classical approach to the study of dynamical systems consists in representing the dynamics of the system in the form of a "source-sink", that means identifying an attractor-repeller pair, which are attractor-repellent sets for all other trajectories of the system. If there is a way to choose this pair so that the space orbits in its complement (the characteristic space of orbits) is connected, this creates prerequisites for finding complete topological invariants of the dynamical system. It is known that such a pair always exists for arbitrary Morse-Smale diffeomorphisms given on any manifolds of dimension $n \geqslant 3$. Whereas for $n=2$ the existence of a connected characteristic space has been proved only for orientation-preserving gradient-like (without heteroclinic points) diffeomorphisms defined on an orientable surface. In the present work, it is constructively shown that the violation of at least one of the above conditions (absence of heteroclinic points, orientability of a surface, orientability of a diffeomorphism) leads to the existence of Morse-Smale diffeomorphisms on surfaces that do not have a connected characteristic space of orbits.
\end{abstract}

\section{Introduction and the statement of the result}
Let $f:M^n\to M^n$ be a Morse-Smale diffeomorphism defined on a closed connected $n$-manifold. Denote by $\Omega^0_f,\,\Omega^1_f,\,\Omega^2_f$ the set of sinks, saddles, and sources of the diffeomorphism $f$. For any (possibly empty) $f$-invariant set $\Sigma\subset\Omega^1_f$ such that $cl(W^u_\Sigma)\setminus W^u_\Sigma\subset\Omega^0_f$, set $$A_\Sigma=\Omega^0_f\cup W^u_\Sigma,\,R_\Sigma=\Omega^2_f\cup W^s_{\Omega^1_f\setminus\Sigma}.$$
It follows from \cite{GMPZ2010} that $A_\Sigma$ and $R_\Sigma$ are an attractor and a repeller, which are called {\it dual}.
In the monograph \cite{grin} the set $$V_\Sigma=M^n\setminus(A_\Sigma\cup R_\Sigma)$$ is called the {\it characteristic space}, and the orbit space $\hat V_\Sigma$ of the action $f$ on $V_\Sigma$ is called the {\it characteristic space of orbits}.

There are a number of examples where a reasonable choice of a dual pair leads to a complete topological classification of some subset of Morse-Smale dynamical systems (look, for example, \cite{BGP}, \cite{BonattiGrinesMedvedevPecou2004}, \cite{GGP}, \cite{GGPM}, \cite{MMP}, and an overview of \cite{GGPZ}). In most cases, finding complete topological invariants is based on the existence of a connected characteristic space of orbits for the class of systems under consideration. For example, according to \cite{BGP}, for any Morse-Smale 3-diffeomorphism, the characteristic space of orbits constructed for the set $\Sigma$ of saddle points with a one-dimensional unstable manifold is connected. This fact played a key role in obtaining a complete topological classification of such diffeomorphisms, obtained in \cite{BGP}. According to \cite{GMPZ2010}, any Morse-Smale diffeomorphism defined on a manifold of dimension $n>3$ also has a connected characteristic space of orbits. For orientation-preserving gradient-like (without heteroclinic points) diffeomorphisms on surfaces there is a result in the work \cite{NoPoxa} that the existence of a connected characteristic orbit space $\hat V_\Sigma$ homeomorphic to the two-dimensional torus $\mathbb T^2$.

The main result of the work is the proof of the fact that the violation of at least one of the conditions (absence of heteroclinic points, orientability of the surface, orientability of the diffeomorphism) leads to the existence of Morse-Smale diffeomorphisms on the surface that do not have a connected characteristic space of orbits. Exactly, the following theorem is proved.

\begin{theorem}\label{th} $ $
\begin{enumerate}
\item On any orientable surface $M^2$ there exists an orientation-changing gradient-like diffeomorphism that does not have a connected characteristic space of orbits.
\item On any non-orientable surface $M^2$ there exists a gradient-like diffeomorphism that does not have a connected characteristic space of orbits.
\item On any surface $M^2$ there exists a Morse-Smale diffeomorphism with heteroclinic points that does not have a connected characteristic space of orbits.
\end{enumerate}
\end{theorem}

{\bf Acknowledgment}: The work was supported by the Russian Science Foundation, project no. 17-11-01041, except for section 3, supported by Laboratory of Dynamical Systems and Applications NRU HSE, grant of the Ministry of science and higher education of the RF, ag. № 075-15-2022-1101

\section{Required information and facts}
Let $M^n$ be a smooth closed orientable manifold and $f$ a diffeomorphism on $M^n$.
For a diffeomorphism $f$, a point $x\in X$ is called {\it wandering} if there exists an open neighborhood $U_x$ of $x$ such that $f (U_x)\cap U_x=\emptyset$. Otherwise, the point $x$ is called {\it non-wandering}. It is immediate from the definition that any point in the neighborhood $U_x$ of a
wandering point $x$ is wandering itself and therefore the set of wandering points is open while the set of non-wandering points is closed. 

The set of all nonwandering points of the diffeomorphism $f$ is called the {\it non-wandering set} and usually denoted by $\Omega_f$.

The simplest examples of hyperbolic sets are primarily the hyperbolic fixed points of a diffeomorphism, which can be classified as follows.
Let $f:X\to X$ be a diffeomorphism and $f(p)=p$. A point $p$ is {\it hyperbolic} if and only if the absolute value of each eigenvalue of the Jacobi matrix  $\left(\frac{\partial f}{\partial x}\right)\vert_{p}$ is not equal to 1. If the absolute values of all the eigenvalues are less than 1, then $p$ is called a {\it attracting (a sink point, or sink)}; if the absolute
values of all the eigenvalues are greater than 1 then $p$ is called a {\it repelling (a source point, or source)}. Attracting or repelling points are called a {\it nodes}. A hyperbolic fixed point that is not a {\it node} is called a {\it saddle point or saddle}.

If the point $p$ is a periodic point $f$ with period $per(p)$, then applying the previous construction to the diffeomorphism $f^{per(p)}$, we obtain a classification of hyperbolic periodic points similar to the classification of fixed hyperbolic points .

The hyperbolic structure of a periodic point $p$ leads to its {\it stable} $W^s_p=\{x\in M^n:\,\lim\limits_{k\to +\infty}d(f^{ kper(p)}(x),p)\to 0\}$ and {\it unstable} $W^u_p=\{x\in M^n:\,\lim\limits_{k\to +\infty }d(f^{-kper(p)}(x),p)\to 0\}$ diversities that are smooth embeddings of $\mathbb R^{n-q_p}$ and $\mathbb R^{q_p}$ respectively. Here $q_p$ is the number of eigenvalues of the Jacobian matrix $\left(\frac{\partial f^{per(p)}}{\partial x}\right)\vert_{p}$ modulo greater than $1$.

For a hyperbolic fixed or periodic point $p$, the stable or unstable manifold is called the {\it invariant manifold} of this point, the connected component of the set $W^u_p\setminus p$ ($W^s_p\setminus p$) is called unstable (stable) {\it separatrix}.

A closed $f$-invariant set $A\subset M^n$ is called an {\it attractor} of a discrete dynamical system $f$ if it has a compact neighborhood $U_A$ such that $f(U_A)\subset int~U_A$ and $A=\bigcap\limits_{k\geq 0}f^k(U_A)$. The neighborhood of $U_A$ is called captivating or isolating. {\it Repeller} is defined as an attractor for $f^{-1}$. An attractor and a repeller are called {\it dual} if the complement to the exciting neighborhood of the attractor is the exciting neighborhood of the repeller.

A diffeomorphism $f:M^n\to M^n$ is called a {\it Morse-Smale diffeomorphism} if

1) the nonwandering set $\Omega_f$ consists of a finite number of hyperbolic orbits;

2) the manifolds $W^s_p$, $W^u_q$ intersect transversally for any nonwandering points $p$, $q$.

A Morse-Smale diffeomorphism is called {\it gradient-like} if the condition $W^s_{\sigma_1} \cap W^u_{\sigma_2}\neq\emptyset$ for different points $\sigma_1, \sigma_2\in \Omega_f$ implies that $dim\, W^u_{\sigma_1}<dim\, W^u_{\sigma_2}.$ In dimension $n=2$, the set of gradient-like diffeomorphisms coincides with the set of Morse-Smale diffeomorphisms whose saddle separatrices do not intersect.

If $M^n$ is an orientable manifold, then the diffeomorphism $f:M^n\to M^n$ is called {\it orientation-preserving}, if $f$ has a positive Jacobian at least at one point, otherwise the diffeomorphism is called {\it orientation-changing}.

Let $f:M^2\to M^2$ be a gradient-like diffeomorphism defined on a closed surface $M^2$. Let $\omega$ be the sink point of the period $m_\omega$ of the diffeomorphism $f$. According to \cite[Theorem 5.5]{palis_melo}, the diffeomorphism $f^{m_\omega}$ in some neighborhood of the point $\omega$ is topologically conjugate to the linear diffeomorphism of the plane given by the matrix
$\begin{pmatrix} \frac{1}{2}&0\cr
	0&\varsigma_\omega\cdot\frac{1}{2}\cr
\end{pmatrix},$ where $\varsigma_\omega=+1\ (-1)$ if $f^{m_\omega}\vert_{W^s_{\omega}}$ preserves (changes) orientation. We say that the sink $\omega$ has a {\it positive orientation type} if $\varsigma_\omega=+1$ and has a {\it negative orientation type} otherwise.

Denote by $\mathcal O_\omega$ the orbit of the point $\omega$. Let $V_\omega=W^s_{\mathcal O_{\omega}}\setminus\mathcal O_{\omega}$. Denote by $\hat V_{\omega}=V_{\omega}/f$ the orbit space of the action of the group $F=\{f^k, k\in\mathbb Z\}\cong\mathbb Z$ on $V_{ \omega}$ and by $p_{_{\omega}}:V_{\omega}\to\hat V_{\omega}$ the natural projection.

\begin{proposition}[\cite{PoKaGr}, Утверждение 1]\label{sos1} 
	The manifold $\hat V_{\omega}$ is diffeomorphic to a two-dimensional torus if $\varsigma_\omega=+1$ and is diffeomorphic to a Klein bottle if $\varsigma_\omega=-1$. Moreover, $\eta_{_\omega}(\pi_1(\hat V_{\omega}))= m_\omega\mathbb Z$.
\end{proposition}

Similarly denote the orientation type $\varsigma_\alpha$ for the periodic source $\alpha$ of the diffeomorphism $f$, the space of orbits $\hat V_\alpha$ and the projection of the stable separatrix of the saddle point into it.

Let $\sigma$ be a saddle point of the period $m_\sigma$ of the diffeomorphism $f$. According to \cite[Theorem 5.5]{palis_melo}, the diffeomorphism $f^{m_\sigma}$ in some neighborhood of the point $\sigma$ is topologically conjugate to the linear diffeomorphism of the plane given by the matrix
$\begin{pmatrix} \nu_\sigma\cdot\dfrac{1}{2}&0\cr
	0&\lambda_\sigma\cdot 2\cr
\end{pmatrix},$
where $\nu_\sigma=+1\ (-1)$ if $f\vert_{W^s_{p}}$ preserves (changes) orientation; $\lambda_\sigma=+1\ (-1)$ if $f\vert_{W^u_{p}}$  preserves (changes) orientation. A pair $\varsigma_\sigma=(\nu_\sigma, \lambda_\sigma)$ will be called the {\it orientation type} of the saddle point $\sigma$ and denote by $a_{\varsigma_\sigma }:\mathbb R^2 \to\mathbb R^2$ a corresponding linear diffeomorphism. If $\nu_\sigma>0,\, \lambda_\sigma>0$, then the type of orientation will be called {\it positive}, and {\it negative} otherwise.

\begin{figure}[h]
\center{\includegraphics[width=0.5\linewidth]{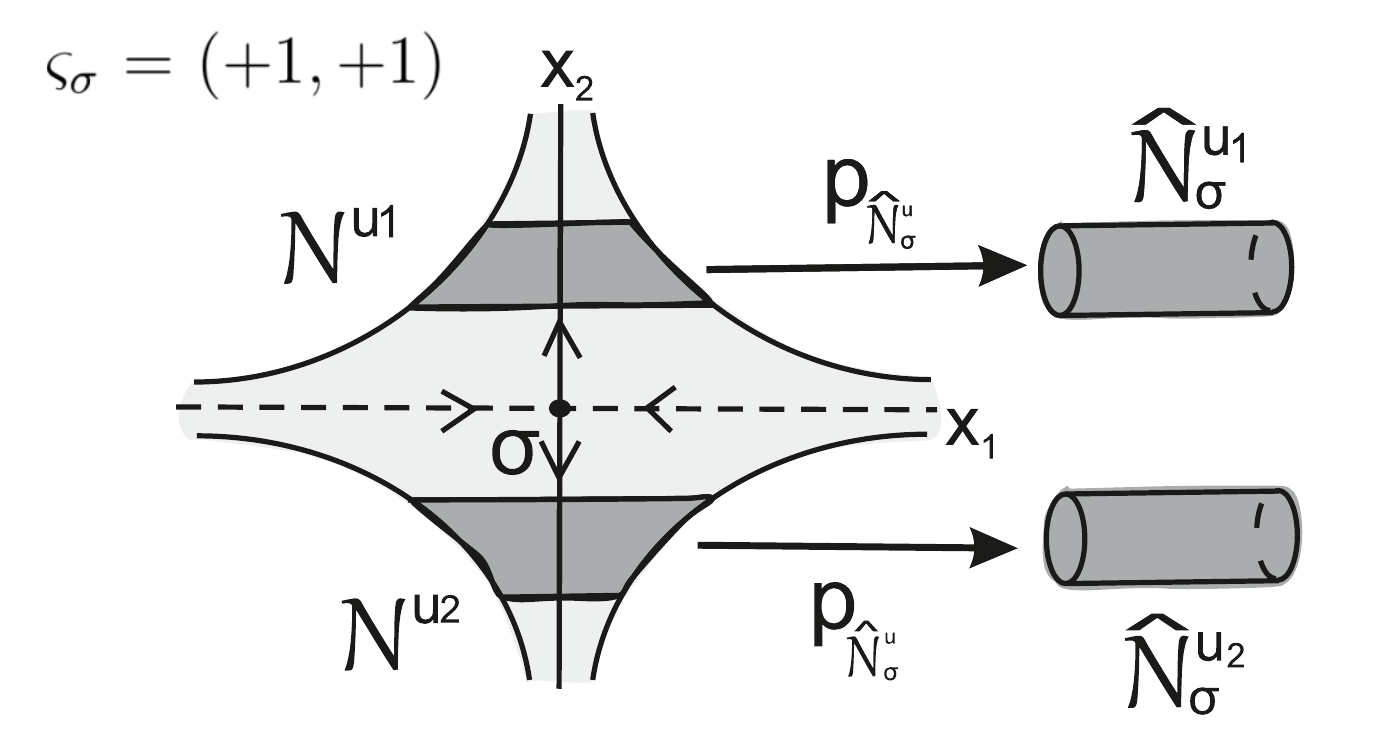}}
	\caption{An orbit space $\hat N^u_{{\sigma}}$}
	\label{perest}
\end{figure}

Denote by ${\mathcal O_\sigma}$ the orbit of the saddle point $\sigma$ and set $N^u_{\sigma}=N_{\mathcal O_\sigma}\setminus W^s_{\mathcal O_\sigma}$ . Then the group $F$ acts on $N^u_{\sigma}$, generating the orbit space $\hat N^u_{{\sigma}}=N^u_{{\sigma}}/f$ and the natural projection $p_ {\sigma}^u:N^u_{{\sigma}}\to\hat N^u_{{\sigma}}$ (see Figure \ref{perest} for the case $\varsigma_{\sigma}=(+1,+1)$).

Denote similarly by the orbit space $\hat N^s_{{\sigma}}=N^s_{{\sigma}}/f$ of the action of the group $F$ on $N^s_{{\sigma}}=N_{\mathcal O_{ \sigma}}\setminus W^u_{\mathcal O_{\sigma}}$, the natural projection
$p_{\sigma}^s:N^s_{{\sigma}}\to\hat N^s_{{\sigma}}$ and mapping $\eta^s_{{\sigma}}$ composed of homomorphisms into the group $\mathbb Z$ from the fundamental group of each connected component of the space $\hat N^s_{{\sigma}}$.

\begin{figure}[h]
\center{\includegraphics[width=0.6\linewidth]{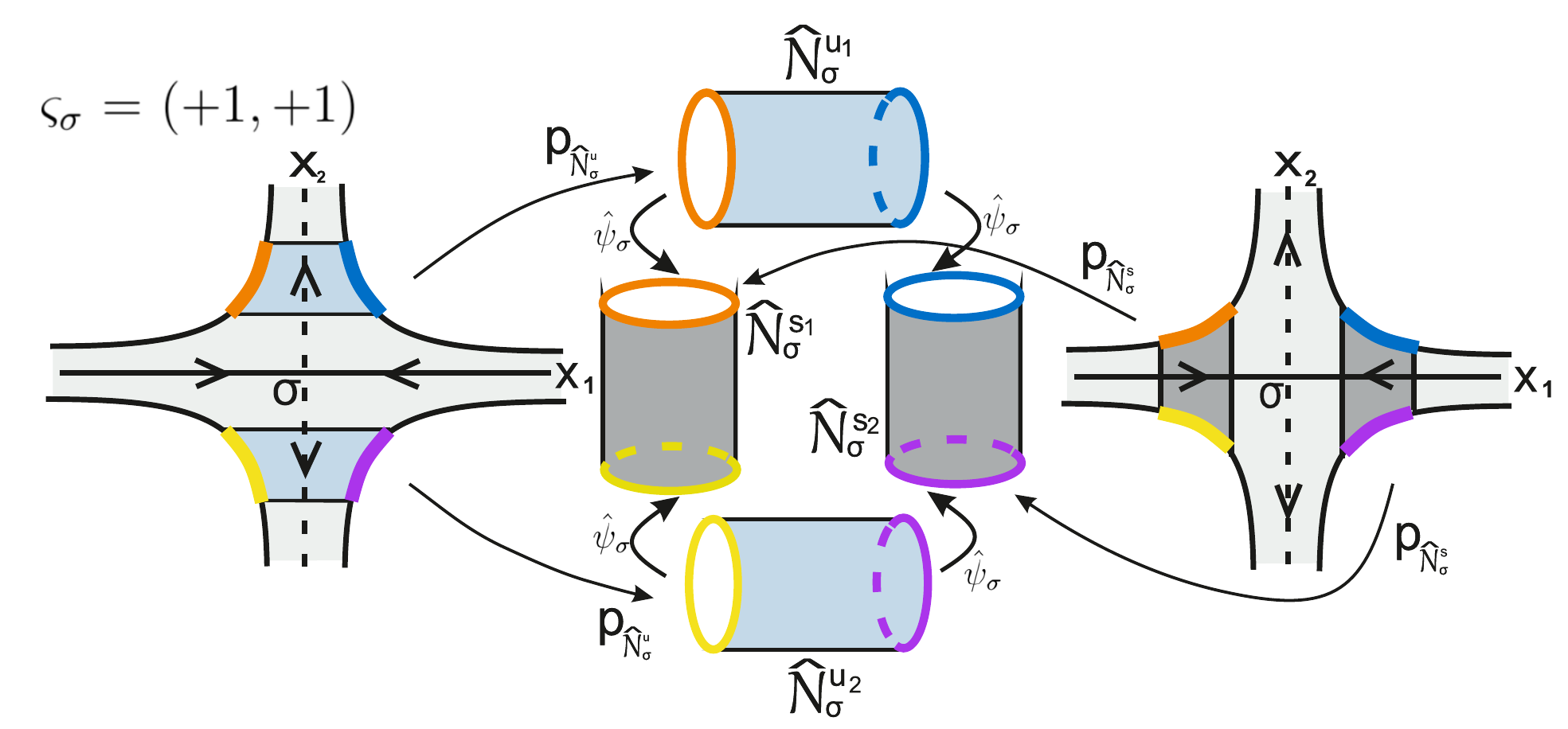}}
	\caption{Rearrangement map}
	\label{operest}
\end{figure}

In addition, the map $$\hat\psi_\sigma=p_{\sigma}^s (p^u_\sigma)^{-1}:\partial {\hat {N}^{u}_{{\sigma}}}\to\partial {\hat {N}^{s}_{{\sigma}}}$$ is well defined and it will be called the {\it rearrangement map} (see Figure \ref{operest} for the case $\varsigma_{\sigma}=(+1,+1)$). 

Denote by $\Omega^0_f,\,\Omega^1_f,\,\Omega^2_f$ the set of sinks, saddles, and sources of the diffeomorphism $f$. For any (possibly empty) $f$-invariant set $\Sigma\subset\Omega^1_f$ such that $cl(W^u_\Sigma)\setminus W^u_\Sigma\subset\Omega^0_f$, we set $$A_\Sigma=\Omega^0_f\cup W^u_\Sigma,\,R_\Sigma=\Omega^2_f\cup W^s_{\Omega^1_f\setminus\Sigma}.$$

The set $$V_\Sigma=M^2\setminus(A_\Sigma\cup R_\Sigma)$$ is called the {\it characteristic space}. The factor space $$\hat V_\Sigma=V_\Sigma/f$$ is called the {\it characteristic space of orbits}. Let $V_\Sigma^1=p_{_\Sigma}^{-1}(\hat V_\Sigma^1),\dots,V_\Sigma^k=p_{_\Sigma}^{-1}( \hat V_\Sigma^k)$ and denote by $m_1,\dots,m_k$ the number of connected components in the sets $V_\Sigma^1,\dots, V_\Sigma^k$ respectively.

\begin{proposition}[\cite{Pix}, Proposition 1]\label{sigma} Each connected component of the characteristic orbit space $\hat V_{\Sigma}$ is homeomorphic either to a two-dimensional torus or to a Klein bottle.
\end{proposition}

\begin{proposition}[\cite{NoPoTsa} Lemma 4.1]\label{L4+}Let $\Sigma'=\Sigma\cup\mathcal O_\sigma$ for some saddle orbit $\mathcal O_\sigma$ and $\hat v,\,\hat v'$ be the disjoint union of the connected components of the spaces $\hat V_{\Sigma},\,\hat V_{\Sigma'}$ which have non-empty intersection with $\hat N^u_\sigma,\,\hat N^s_\sigma$, respectively. Then $$\hat V_{\Sigma'}\cong(\hat V_\Sigma\setminus int\,\hat N^u_\sigma)\cup_{\hat\psi_\sigma}\hat N^s_\sigma. $$ Wherein
\begin{equation}\label{vv}
\hat V_{\Sigma'}\cong(\hat V_\Sigma\setminus \hat v)\sqcup\hat v'
\end{equation}
\end{proposition}

\begin{figure}[h]
\center{\includegraphics[width=1\linewidth]{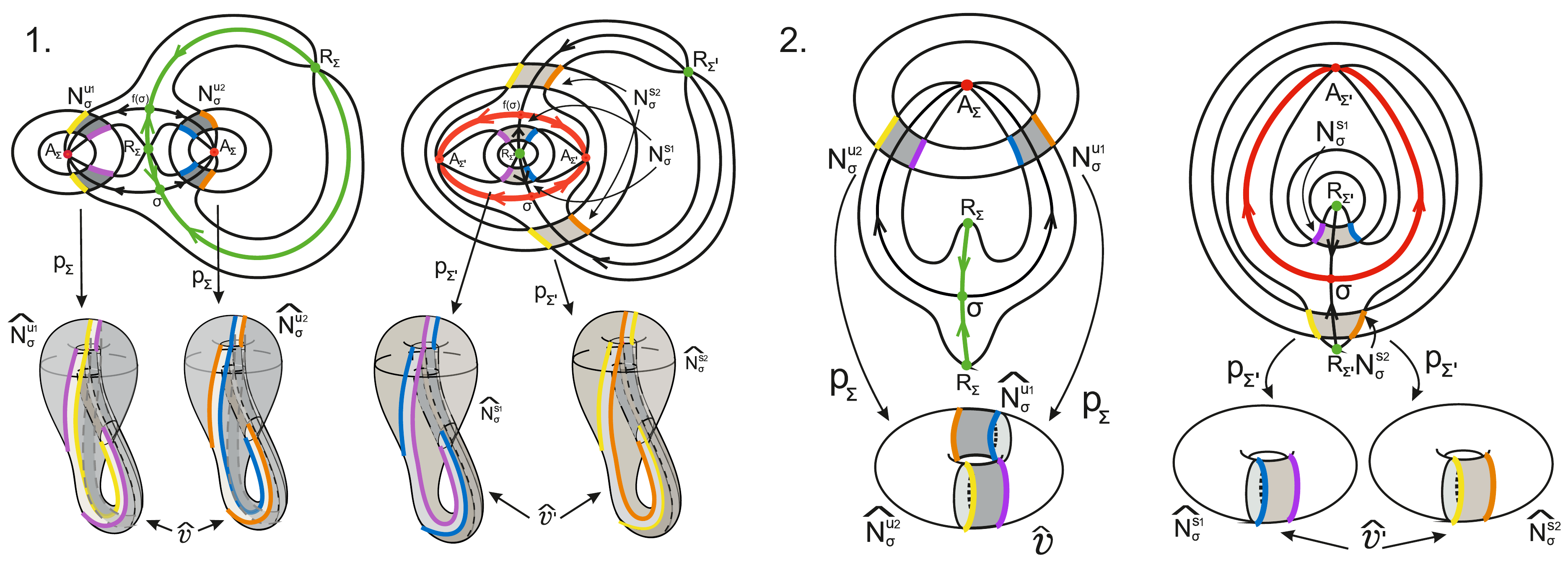}}
	\caption{Illustration for the collary \ref{L4}}
	\label{LemmaP}
\end{figure}

\begin{collary}\label{L4} If $\sigma$ is a saddle point with a positive orientation $\varsigma_{\sigma}=(+1,+1)$, then for the sets $\hat v,\,\hat v'$ in the formula (\ref {vv}) the following features are implemented:
\begin{itemize}
\item $\hat v$ -- a disjoint union of two Klein bottles, $\hat v$ -- a disjoint union of two Klein bottles (see Figure \ref{LemmaP}(1));
\item $\hat v$ -- a torus, $\hat v'$ -- a disjoint union of two tori (see Figure \ref{LemmaP}(2));
\item $\hat v$ -- a disjoint union of two tori and $\hat v'$ -- a torus;
\item $\hat v$ -- a disjoint union of a torus and a Klein bottle and $\hat v'$ — a Klein bottle;
\item $\hat v$ -- a Klein bottle and $\hat v'$ -- a disjoint union of a torus and a Klein bottle;
\item $\hat v$ -- a torus and $\hat v'$ — a torus (if $M^2$-nonorientable surface).
	\end{itemize}
If $\sigma$ is a saddle point with negative orientation $\varsigma_{\sigma}=(-1,-1)$, then the following possibilities are realized:
\begin{itemize}
	\item $\hat v$ — a Klein bottle and $\hat v'$ — a Klein bottle;
	\item $\hat v$ — a torus and $\hat v'$ -- a torus.
\end{itemize}
\end{collary}

\section{Construction of model diffeomorphisms}
In this section, we construct several basic diffeomorphisms, the proof of the theorem \ref{th} will be based on them.

\subsection{Gradient-like diffeomorphism $\psi_0$ on the sphere $\mathbb S^2$}\label{phi}

\begin{figure}[h]
	\center{\includegraphics[width=0.3\linewidth]{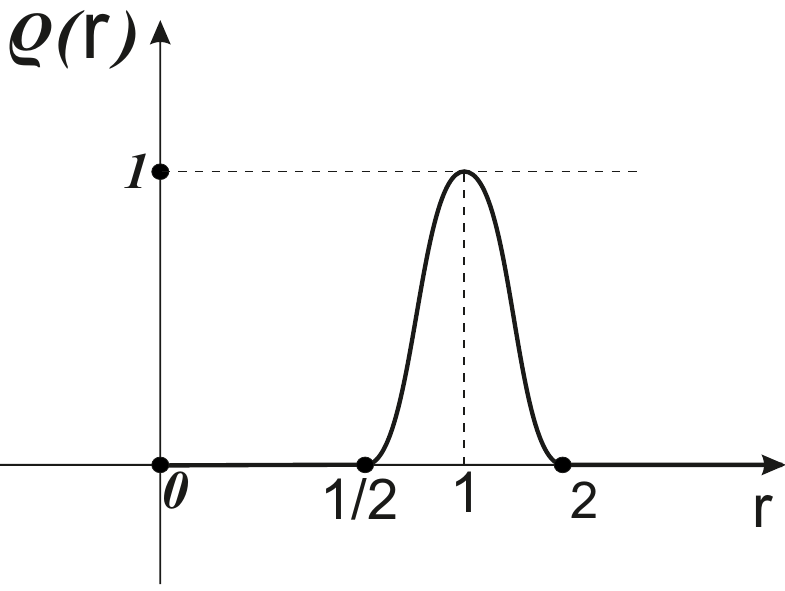}}
	\caption{Function graph $\varrho$.}
	\label{ris:varrhofnuc}
\end{figure}

\begin{figure}[h]
	\center{\includegraphics[width=0.4\linewidth]{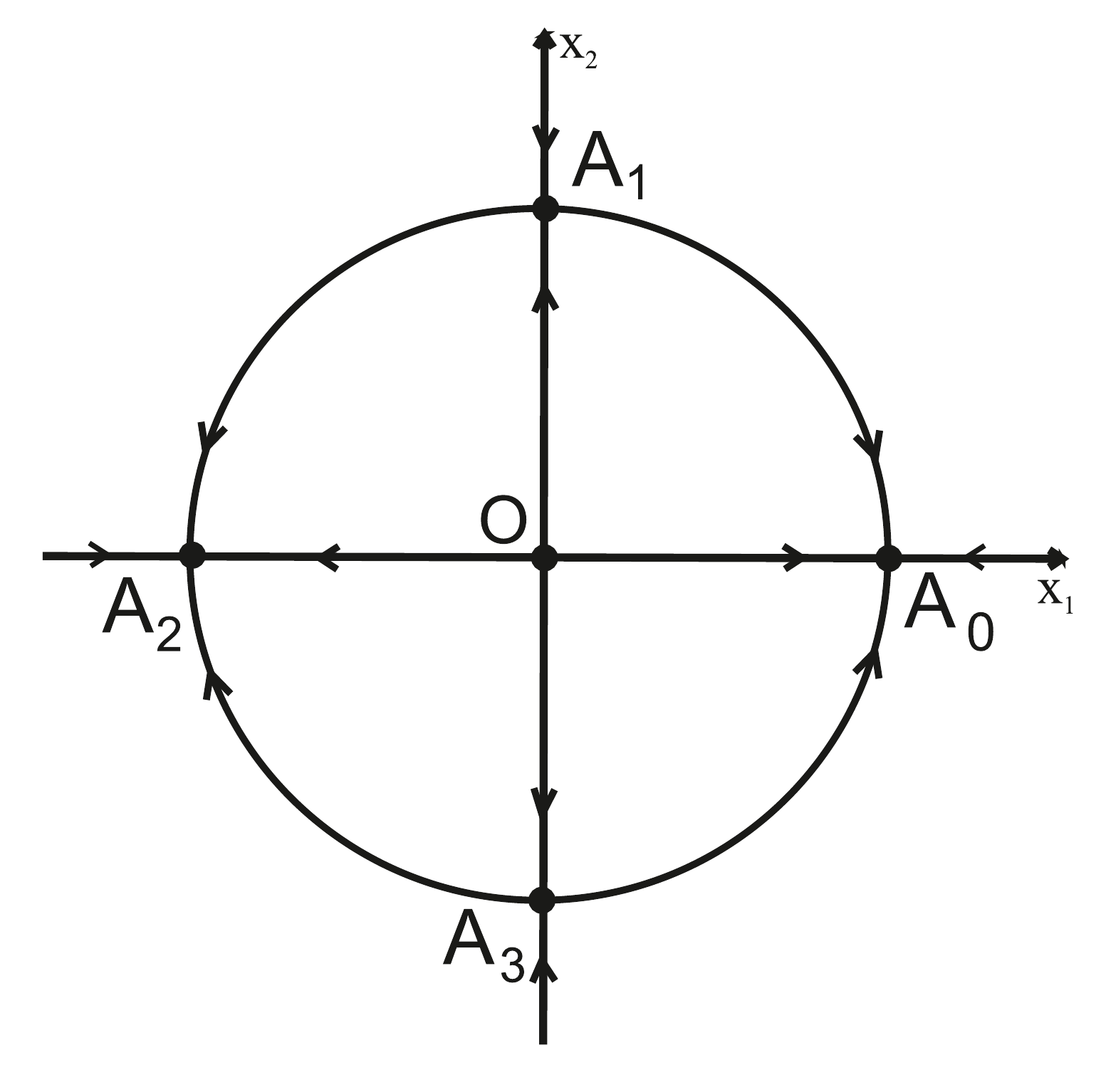}}
	\caption{Phase portrait of a diffeomorphism $\chi$}
	\label{realpic}
\end{figure}

Define polar coordinates $(r,\varphi)$ on the plane $\mathbb R^2$.
Denote by $\varrho(r)$ the function depicted on the graph(see Figure \ref{ris:varrhofnuc}),which has the property $\varrho(r)=\varrho(\frac 1r)$.
Also define a vector field on the plane $\mathbb R^2$ using the following system of differential equations:
$$\begin{cases}\dot r=\begin{cases}-r(r-1), 0\leqslant r \leqslant 1;\\ 1-r, r>1;\end{cases}\\ \dot\varphi= \varrho(r)sin 2\varphi.\end{cases}$$ 
Denote by  $\chi^t$ the flow induced by this vector field, and denote by  $\chi^t$ the diffeomorphism, which is the shift of the flow $\chi^t$ per unit of time. The result is a diffeomorphism that has a hyperbolic source at the origin $O$, hyperbolic saddles at points $A_1, A_3$ and hyperbolic drains at points $A_0, A_2$ (see Figure \ref{realpic}) .

\begin{figure}
\centering{\includegraphics[width=8cm,height=4.5 cm]{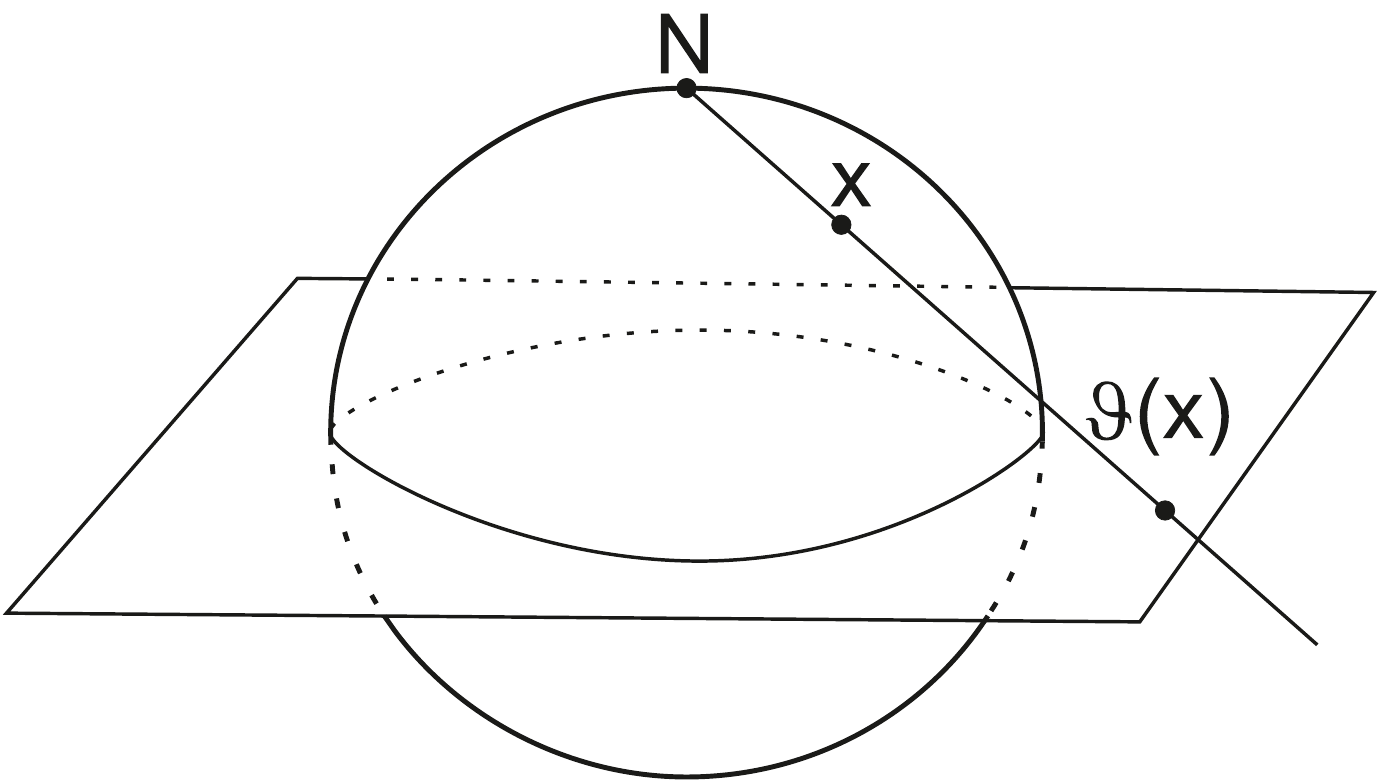}} 
\caption{Stereographic projection}\label{stereo}
\end{figure}

Let a diffeomorphism $\theta: \mathbb R^2 \to \mathbb R^2$ be as follows as $\theta (r,\varphi)=(r, -\varphi)$.
We define the diffeomorphism $\bar \psi_0:\mathbb R^2 \to \mathbb R^2$ by the formula $$\bar \psi_0=\theta\circ\chi.$$ By the construction, the nonwandering set of the diffeomorphism $\bar \psi_0$ coincides with a nonwandering diffeomorphism set $\chi$.
Consider the standard two-dimensional sphere $$\mathbb S^2=\{(x_1,x_2,x_3)\in\mathbb R^3:x_1^2+x^2_2+x_3^2=1\}.$$ Denote by $ N(0,0,1)$ north pole and define a stereographic projection (see Figure \ref{stereo}) $\vartheta:\mathbb {S}^2\setminus\{N\}\to\mathbb{R }^2$ formula
$$\vartheta(x_1,x_2,x_3)=\left(\frac{x_1}{1-x_3},\frac{x_2}{1-x_3}\right).$$

\begin{figure}[h]
\centerline{\includegraphics
		[width=7 cm]{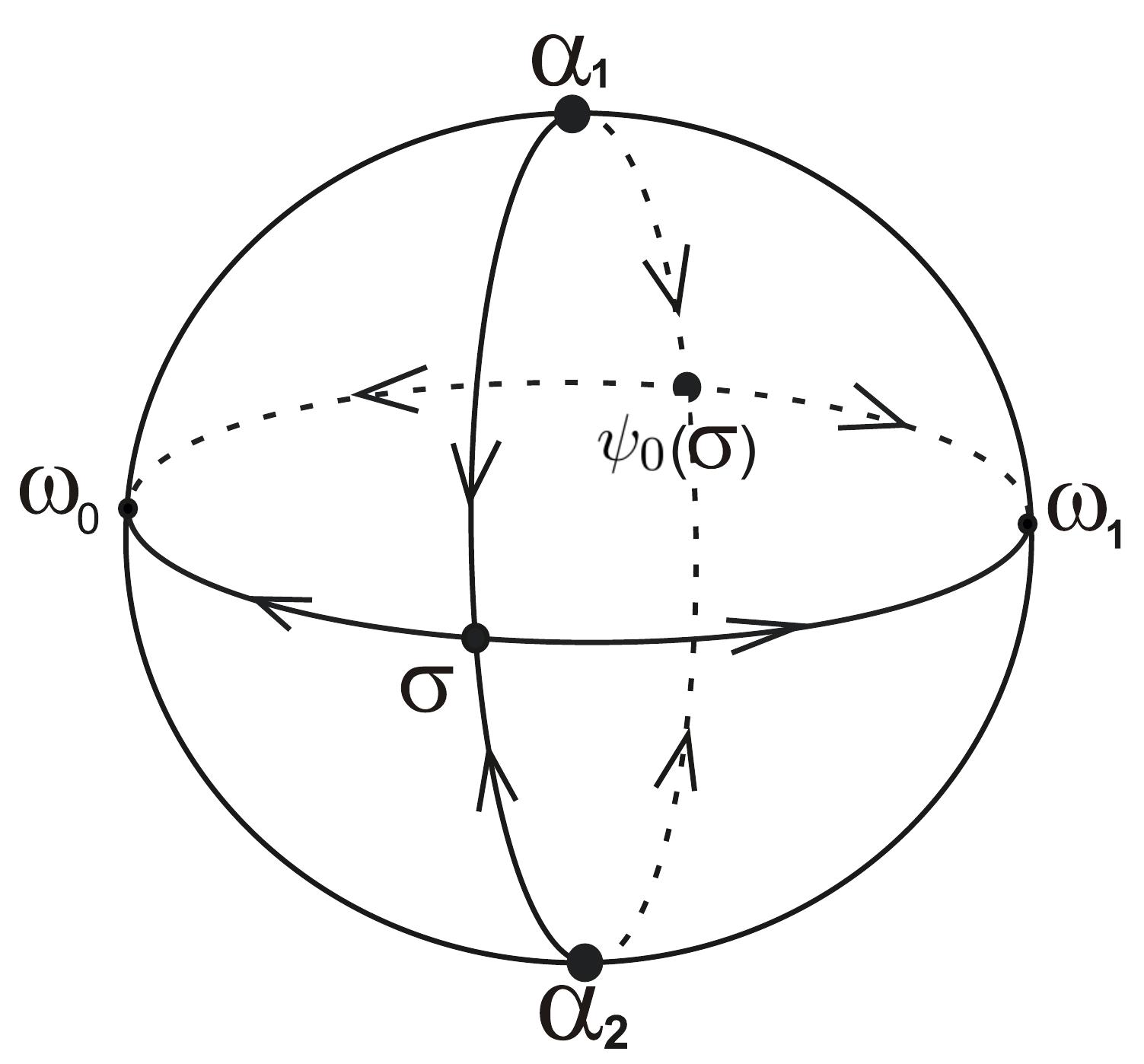}}
	\caption{\small Diffeomorphism $\psi_0$}\label{SAR}
\end{figure}

Define a diffeomorphism $\psi_0: \mathbb S^2 \to\mathbb S^2$ by the formula $$\psi_0(x)=\begin{cases}
	\vartheta^{-1}\circ \bar \psi_0 \circ \vartheta(x),\,x\in S^2\setminus\{N\}, \\
	N,\,x=N.
\end{cases}$$
By construction, $\psi_0$ is an orientation-changing gradient-like 2-sphere diffeomorphism whose nonwandering set consists of two fixed sources $\alpha_1=N$, $\alpha_2=\vartheta^{-1}(O)$ of negative orientation $(\varsigma_{\alpha_1}=\varsigma_{\alpha_2}=-1)$; two fixed sinks $\omega_0=\vartheta^{-1}(A_0)$, $ \omega_1=\vartheta^{-1}(A_{2})$ of negative orientation $(\varsigma_{\omega_0}=\varsigma_{\omega_1}=-1)$ and saddle orbit $\mathcal O_\sigma=\{\sigma=\vartheta^{-1}(A_1),\,\psi_0(\sigma)=\vartheta^{-1}(A_{3})\}$ of period $2$ with orientation type $\varsigma_\sigma=(+1,+1)$  (see Figure \ref{SAR}): $$\Omega_{\psi_0}=\{\alpha_1,\alpha_2,\omega_0,\omega_1,\sigma,\psi_0(\sigma)\}.$$

\subsection{Gradient-like diffeomorphism $\tilde\psi_1$ on the projective plane $\mathbb R P^2$}\label{psi0}

\begin{figure}[h]
\center{\includegraphics[scale=0.35]{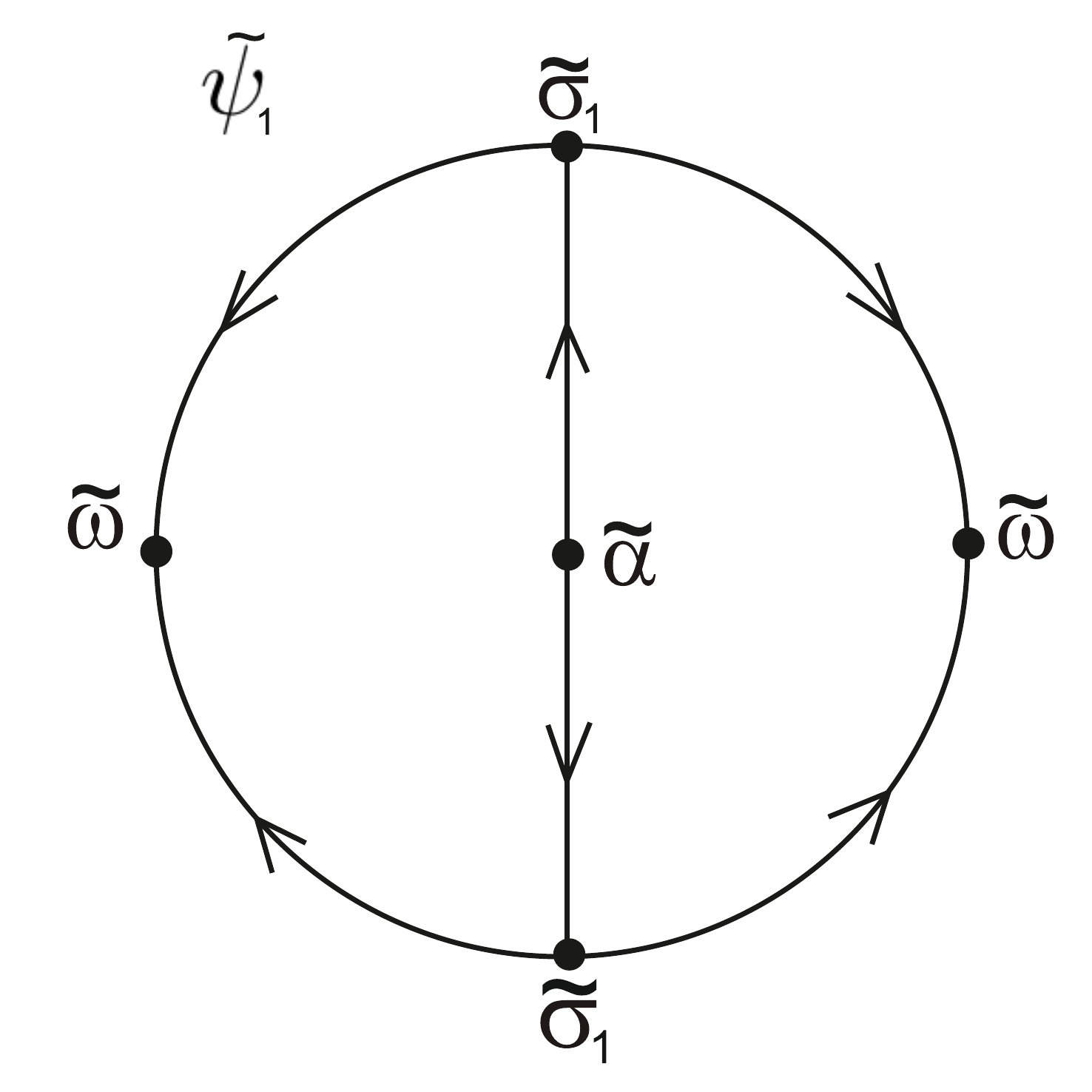}}
	\caption{\small Diffeomorphism $\tilde\psi_1$}
	\label{Diffpsi}
\end{figure}

Consider the diffeomorphism $\psi_0:\mathbb S^2 \to \mathbb S^2$ defined in \ref{phi} and the group $\mathbb Z_2=\{ +1, -1\} $ acting on the two-dimensional sphere $\mathbb S^2=\{(x_1,x_2,x_3)\in\mathbb R^3:x_1^2+x^2_2+x_3^2=1\}$ as follows:
$$\pm{1} \cdot x=\pm{x},\,  x=(x_1, x_2, x_3)\in \mathbb S^2 .$$
Then the orbit space $\mathbb S^2/ \mathbb Z_2$ of the action of the group $\mathbb Z_2$ on $\mathbb S^2$ is the projective plane $\mathbb R P^2$. Let $p:\mathbb S^2 \rightarrow \mathbb R P^2$ be the natural projection. Let the diffeomorphism $\tilde\psi_1: \mathbb R P^2 \rightarrow \mathbb R P^2$ be defined by the formula $$\tilde\psi_1(x)=p \circ \psi_0 \circ p^{-1}(x),\ , x\in \mathbb R P^2.$$ By construction, the nonwandering set of the constructed diffeomorphism consists of three fixed points: the source $\tilde\alpha$ of negative orientation $(\varsigma_{\tilde\alpha}=-1)$, the sink $\tilde\omega$ of negative orientation $(\varsigma_{\tilde\omega}=-1)$ and saddle $\tilde\sigma_1$ with orientation type $\varsigma_{\tilde\sigma_1}=(-1,-1)$ (see Figure \ref{Diffpsi}): $$\Omega_{\tilde\psi_1}=\{\tilde \alpha, \tilde \omega, \tilde \sigma_1\}.$$

\subsection{Gradient-like diffeomorphism $\tilde\psi_q$ on a non-orientable surface of genus q}

\begin{figure}[h]
\center{\includegraphics[scale=0.35]{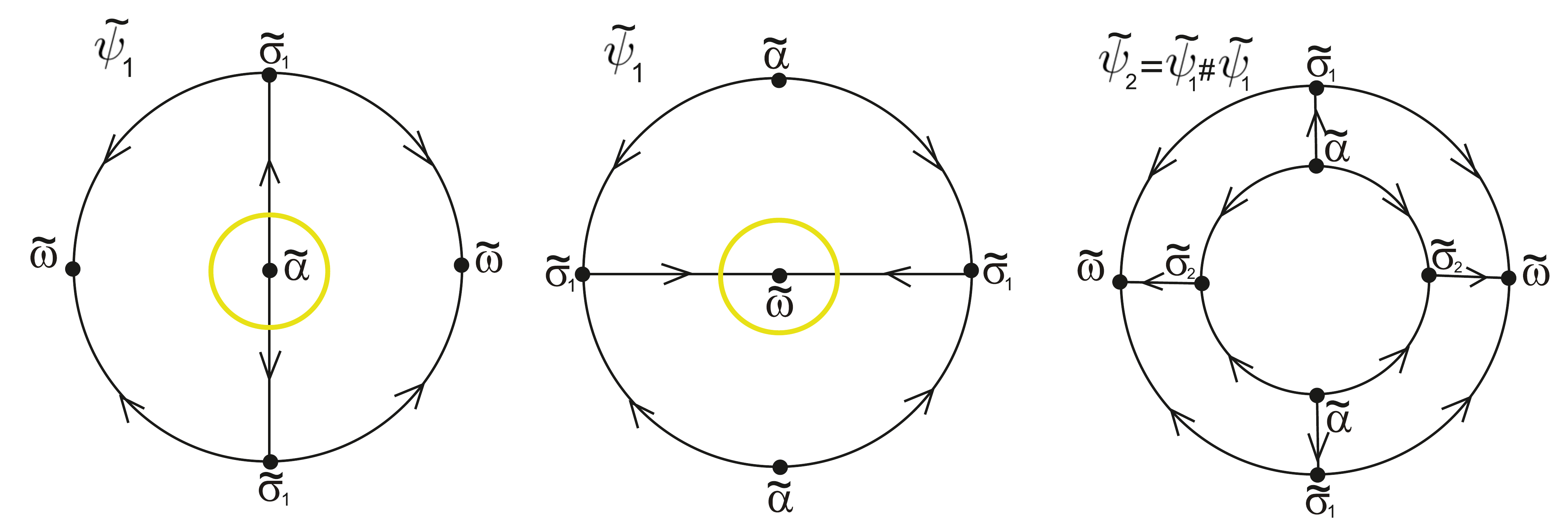}}
\caption{\small Diffeomorphism $\tilde\psi_2$}
\label{Diffpsi_q}
\end{figure}

Let $S^-_q=\mathbb S^2\sharp\underbrace{\mathbb RP^2\sharp\dots\sharp\mathbb RP^2}_q$. Construct a model diffeomorphism $\tilde\psi_q:S^-_q\to S^-_q$. The source $\tilde \alpha$ and the sink $\tilde \omega$ of the diffeomorphism $\tilde\psi_1$ should be considered. Whereas they are hyperbolic, there are non-intersecting 2-discs around them $B_{\tilde \omega},\,B_{\tilde \alpha}$, что $\tilde \psi_1 (B_{\tilde \omega}) \subset int\,B_{\tilde \omega},\, \tilde \psi_1 ^{-1}(B_{\tilde \alpha})\subset int\, B_{\tilde \alpha}$. 
Then the connected sum of two copies of projective planes along the disks $B_{\tilde \omega},\,B_{\tilde \alpha}$ is a non-orientable surface $S^-_2$ of genus 2 (see Figure \ref{Diffpsi_q}). Since the dynamics in the disk $B_{\tilde \omega}$ is inverse to the dynamics in the disk $B_{\tilde \alpha}$, a diffeomorphism $\tilde\psi_2$ is well defined on the surface $S^-_2$, which coincides with $\tilde\psi_1$ to ${\mathbb{RP}}^2 \setminus B_{\tilde \alpha}$ and to ${\mathbb{RP}}^2 \setminus B_{\tilde \omega}$. We say that the diffeomorphism $\tilde\psi_2$ is the {\it connected sum} of two copies of the diffeomorphism $\tilde\psi_1$ ($\tilde\psi_2=\tilde\psi_1\sharp\tilde\psi_1$) along the sink $\tilde\omega$ and source $\tilde\alpha$. 

By induction, a diffeomorphism $\tilde\psi_q:S^-_q \rightarrow S^-_q$ on a nonorientable surface of genus $q \geqslant 2$ is constructed as a connected sum of diffeomorphisms $\tilde\psi_{q-1}$ and $\tilde \psi_1$ ($\tilde\psi_q=\tilde\psi_{q-1}\sharp \tilde\psi_1$) along the $\tilde\omega$ sink and $\tilde\alpha$ source. By construction, the nonwandering set of the diffeomorphism $\tilde\psi_q$ consists of $q+2$ fixed points: the source $\tilde\alpha$ of negative orientation $(\varsigma_{\tilde\alpha}=-1)$, the sink $\tilde\omega$ negative orientation $(\varsigma_{\tilde\omega}=-1)$ and $q$ saddles $\tilde\sigma_1,\dots,\tilde\sigma_q$ with orientation type $\varsigma_{\tilde\sigma_i}=(+1,+1)$ (see Figure \ref{Diffpsi_q} for $q=2$):  $$\Omega_{\tilde\psi_q}=\{\tilde \alpha, \tilde \omega, \tilde \sigma_1,\dots,\tilde \sigma_q\}.$$

\subsection{Gradient-like diffeomorphism $\psi_1$ on a torus $\mathbb T^2$}\label{f0T}

\begin{figure}[h!]
\centerline{\includegraphics[width= 5 cm]{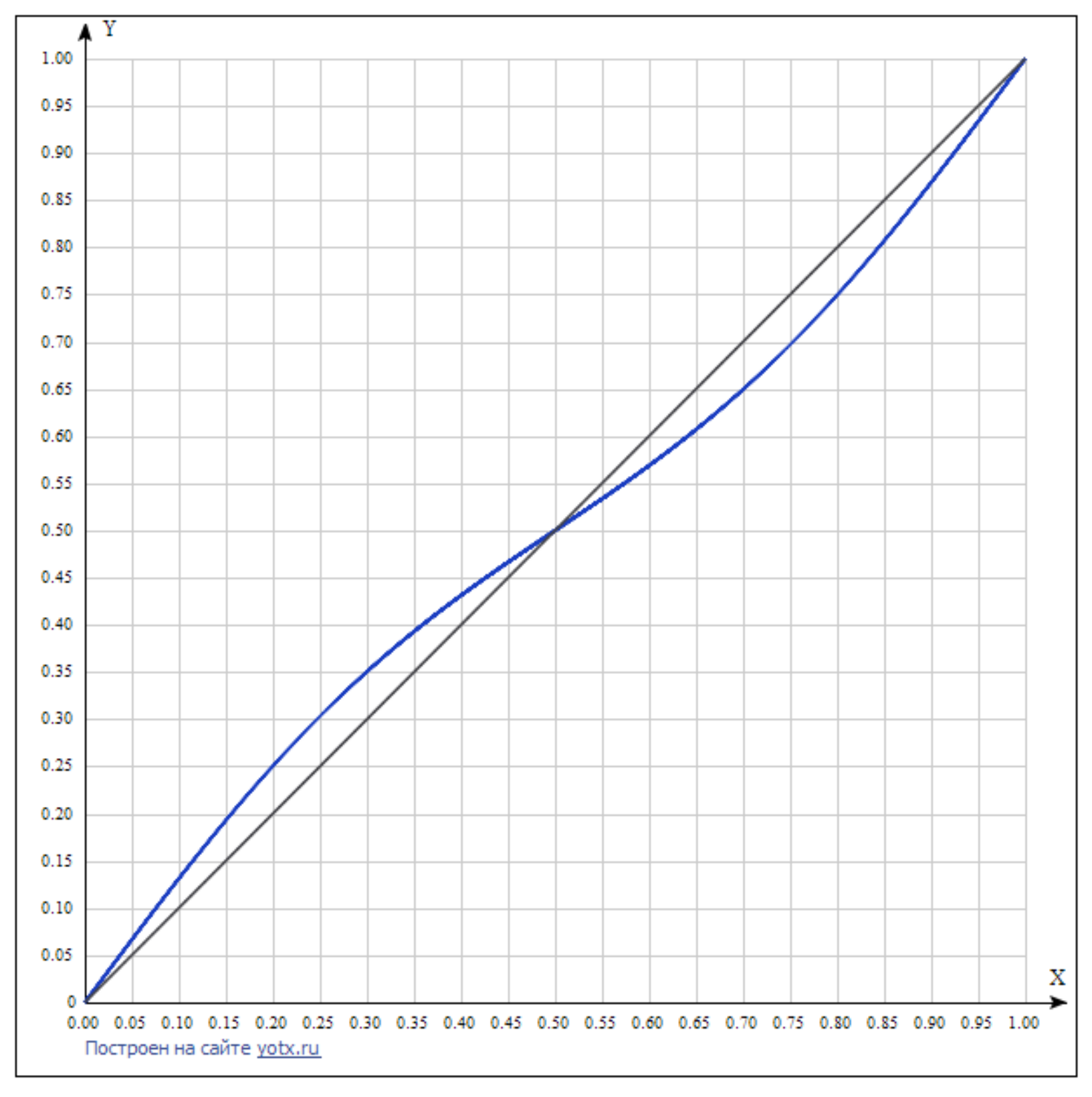}}
	\caption{\small Function Graph $\bar F$}\label{barf0}
\end{figure}

We construct a diffeomorphism $\psi_1$ on the two-dimensional torus $\mathbb T^2$ as a cartesian product of two orientation-preserving source-sink diffeomorphisms on the circle $\mathbb S^1$. For this we should consider the function $\bar{F}:\mathbb R\to\mathbb R$ given by the formula:
$$\bar{F}(x)=x+\frac{1}{6\pi} sin 2\pi x$$ (see Figure \ref{barf0}).

Consider the projection $\pi: \mathbb{R}\to \mathbb{S}^1$ given by the formula $\pi(x)=e^{2\pi i x}$. Since the function $\bar{F}$ is strictly monotonically increasing and satisfies the condition $\bar{F}(x+1)=\bar{F}(x)+1$, it admits a projection onto a circle in diffeomorphism $F:\mathbb S^1\to\mathbb S^1$ given by $$F(z)=\pi\bar F \pi^{-1}(z),\,z\in\mathbb S^1.$$

\begin{figure}[h!]
\centerline{\includegraphics
		[width= 4 cm]{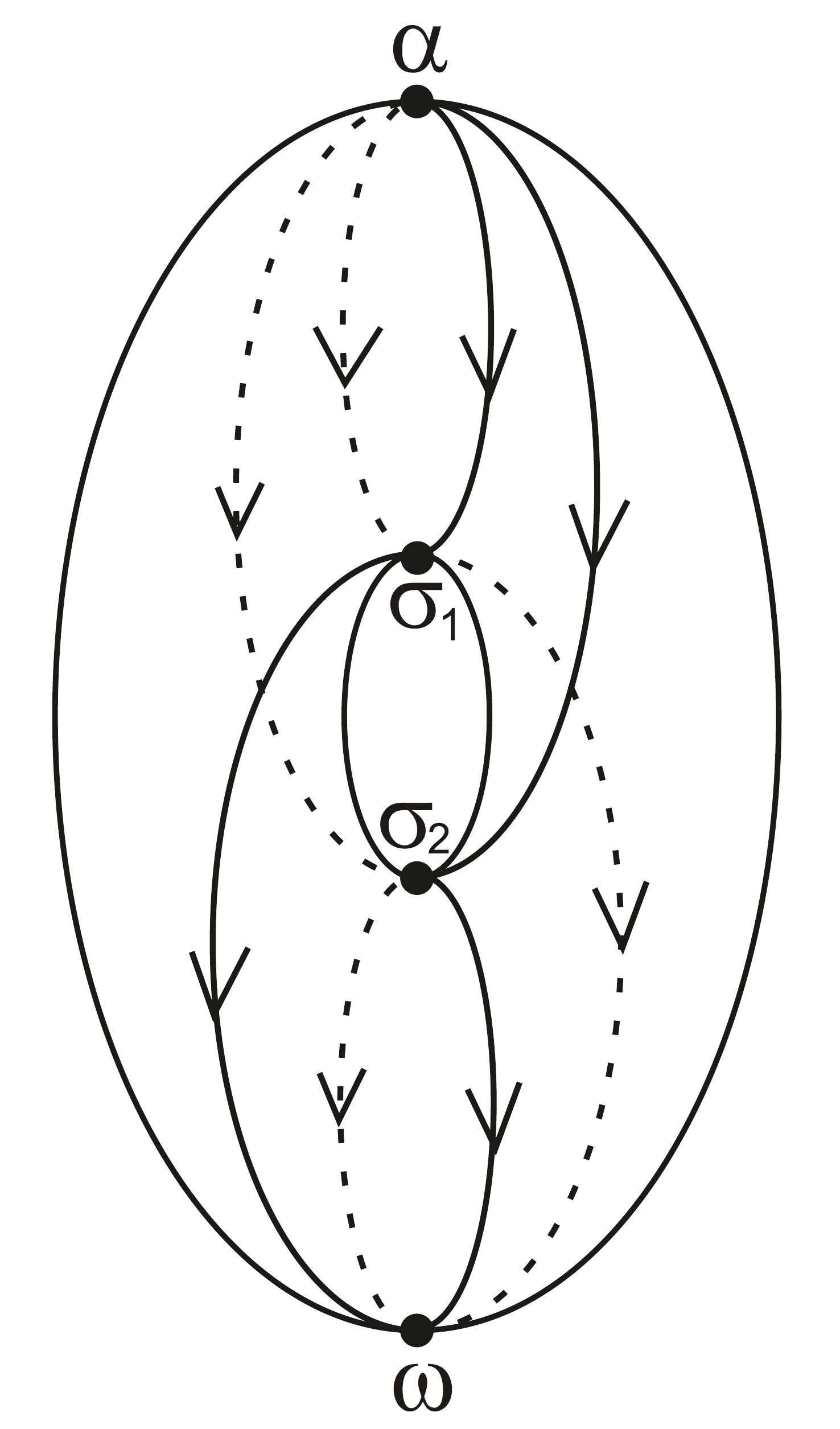}}
	\caption{\small Diffeomorphism $F_1$}\label{f1}
\end{figure}

By construction, the diffeomorphism $F$ has a fixed hyperbolic sink and source and is an orientation-preserving source-sink diffeomorphism. Define a diffeomorphism $F_1:\mathbb T^2\to\mathbb T^2$ by the formula $$F_1(z,w)=(F(z),F(w)),\,z,w\in\mathbb S ^1.$$
Then the diffeomorphism $F_1$ is orientation-preserving, and its nonwandering set consists of four fixed points: a source $\alpha$ of positive orientation $(\varsigma_{\alpha}=+1)$, a sink $\omega$ of positive orientation $( \varsigma_{\omega}=+1)$ and two saddles $\sigma_1,\,\sigma_2$ of positive orientation type (see Figure \ref{f1}): $$\Omega_{F}=\{ \alpha, \omega, \sigma_1, \sigma_2 \}.$$

\begin{figure}[h!]
\centerline{\includegraphics
		[width= 3.5 cm]{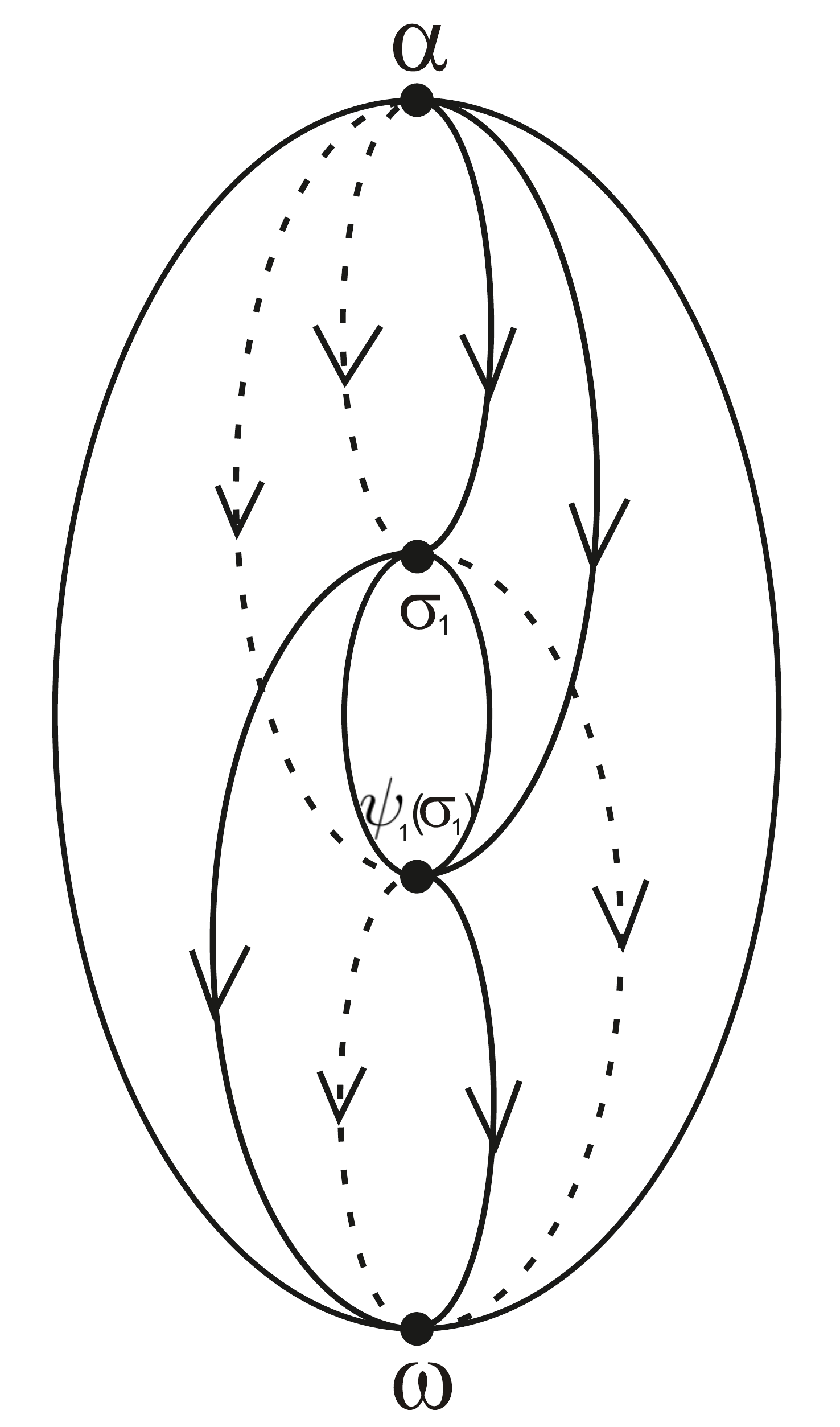}}
	\caption{\small  Diffeomorphism $\psi_1$}\label{psi1}
\end{figure}

Let us represent the two-dimensional torus $\mathbb{T}^2$ as the factor group of the group $\mathbb{R}^2$ with respect to the integer lattice $\mathbb{Z}^2:\mathbb{T}^2=\mathbb{R }^2/\mathbb{Z}^2$. Consider the matrix $A=
\begin{pmatrix}
	0& 1\\
	1& 0
\end{pmatrix}\in GL(2,\mathbb{Z})$ and the algebraic torus automorphism $\hat A:\mathbb T^2\to\mathbb T^2$,
$$\hat A(x,y)=(y,x)\pmod{1}$$
Let $$\psi_1=\hat A\circ F_1: \mathbb{T}^2\to\mathbb T^2.$$
By construction, the diffeomorphism $\psi_1$ is an orientation-changing gradient-like diffeomorphism whose nonwandering set consists of a source $\alpha$ and a sink $\omega$ of negative orientation types $(\varsigma_{\alpha}=\varsigma_{\omega} =-1)$, as well as periodic saddle orbit $\mathcal O_{\sigma_1}=\{\sigma_1,\psi_1(\sigma_1)\}$ of period 2 and orientation type $\varsigma_{\sigma_1}=(+1,+1)$ (see Figure \ref{psi1}):
$$\Omega_{\psi_1}=\{\alpha,\omega,\sigma_1,\psi_1(\sigma_1)\}.$$

\subsection{Gradient-like diffeomorphism $\psi_g$ on an orientable surface of genus $g$}

\begin{figure}[h!]
\centerline{\includegraphics
		[width=12 cm]{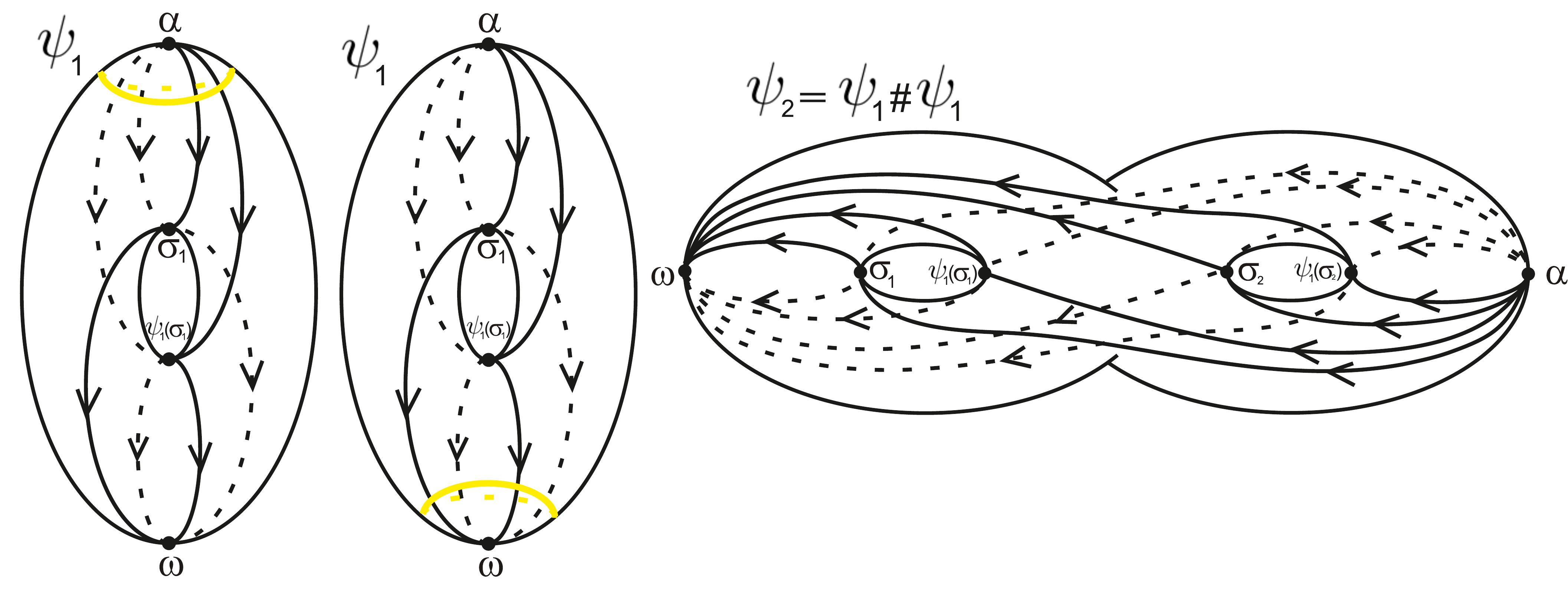}}
	\caption{\small Diffeomorphism  $\psi_2$}\label{SS}
\end{figure}

Let $S^+_g=\mathbb S^2\sharp\underbrace{\mathbb T^2\sharp\dots\sharp\mathbb T^2}_g$. Let us construct a model diffeomorphism $\psi_g:S^+_g\to S^+_g$. To do this, firstly, we should construct a diffeomorphism $\psi_2:S^+_2\to S^+_2$
as a connected sum of two copies of the diffeomorphism $\psi_{1}$ ($\psi_2=\psi_1\sharp \psi_1$) along the sink $\omega$ and the source $\alpha$ (see Fig. \ref{SS}). Then, by induction, we define the diffeomorphism $\psi_g:S^+_g\to S^+_g$
as a connected sum of diffeomorphisms $\psi_{g-1}$ and $\psi_{1}$ ($\psi_{g}=\psi_{g-1}\sharp \psi_1$) along the sink $\omega$ and the source $\alpha$.
By construction, the nonwandering set of the diffeomorphism $\psi_g$ consists of a fixed source $\alpha$ and a fixed sink $\omega$ of negative orientation types $(\varsigma_{\alpha}=\varsigma_{\omega}=-1)$, and $ g$ saddle periodic orbits $\mathcal O_{\sigma_1}=\{\sigma_1,\psi_1(\sigma_1)\},\dots,\mathcal O_{\sigma_g}=\{\sigma_g,\psi_1(\sigma_g) \}$ of period 2 and orientation type $\varsigma_{\sigma_i}=(+1,+1)$ (see Figure \ref{SS} for $g=2$):
$$\Omega_{\psi_g}=\{ \alpha, \omega, \sigma_1,\psi_g(\sigma_1),\dots, \sigma_{g},\psi_g(\sigma_g)\}.$$ 

\subsection{Non-gradient-like Morse-Smale diffeomorphism $\xi_0$ on the sphere $\mathbb S^2$}\label{psi0}

Consider the diffeomorphism $h: \mathbb R^2 \to \mathbb R^2$given in polar coordinates $(r,\varphi),\,\varphi\in\left(-\frac{\pi}{2},\frac{3\pi}{2}\right]$ by the formula $h(r,\varphi)=\left(\frac{r}{2}, \varphi\right)$. Let $A_-=\{(r,\varphi)\in\mathbb R^2\setminus \{0\}:\,\frac{3\pi}{4}\leqslant\varphi\leqslant \frac{5\pi}{4}\}$, $A_+=\{(r,\varphi)\in\mathbb R^2\setminus \{0\}:\,|\varphi|\leqslant \frac{\pi}{4}\}$. Let $C=\mathbb{R}\times[-2,2]$ and define the diffeomorphisms $\eta_-:A_-\to C,\,\eta_+:A_+\to C$ by the formulas $$\eta_-(r,\varphi)=\left(3-log_2r, 8\left(\frac{\varphi}{\pi}-1\right)\right),\,\eta_+(r,\varphi)=\left(-3-log_2r, \frac{8\varphi}{\pi}\right).$$It is directly verified that the diffeomorphism $\eta_-\,(\eta_+)$ conjugates the diffeomorphism $h$ with the diffeomorphism $g:C\to C$ defined by the formula $g(x,d)=(x+1,d) $ and $\eta_-(2,\pi)=(2,0)\,(\eta_+(1/2,0)=(-2,0))$. It is obvious that the diffeomorphism $g$ is included in the flow $g^t:C\to C$ defined by the formula $$g^t(x,d)=(x+t,d).$$ 
{We define the flow $\phi_-^t$ on $C$ using the formulas
	$$\begin{cases}
		\dot{x}=\begin{cases}1-\frac{1}{9}(x^2+d^2-4)^2, \quad x^2+d^2 \leq 4 \cr
			1, \quad\mbox{else}
		\end{cases}\cr
		\dot{d}=\begin{cases}
			\frac{d}{2}\big(\sin\big(\frac{\pi}{2}\big(x^2+d^2-3\big)\big)-1\big), \quad 2<x^2+d^2\leq 4\cr
			-d,\quad \quad x^2+d^2\leq 2\cr
			0, \quad\mbox{else}
		\end{cases}\cr
	\end{cases}$$} 
\begin{figure}[h]
\centerline{\includegraphics[width=9 true cm, height=7  true cm]{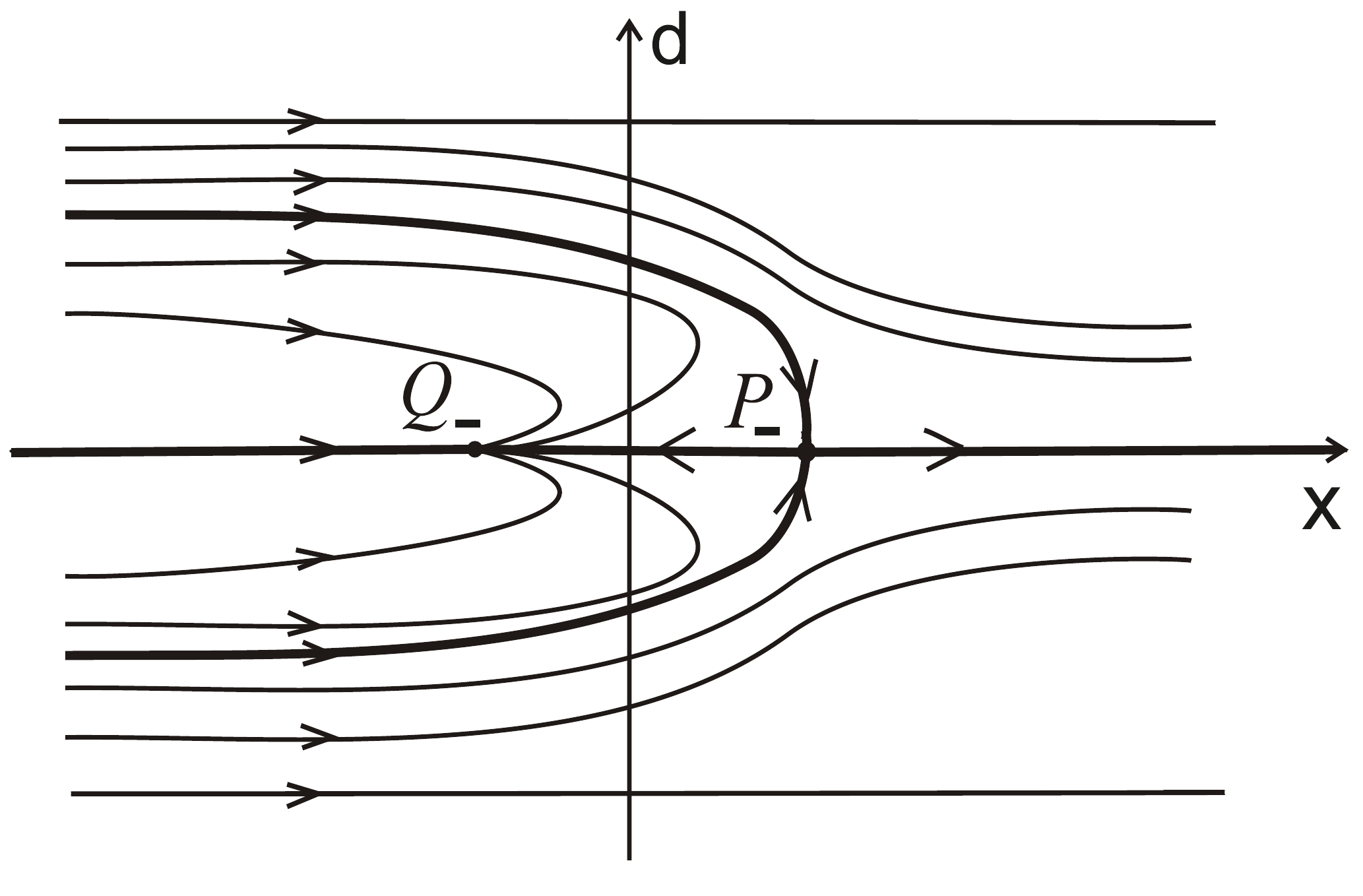}}
\caption{Flow paths $\phi_-^t$}\label{cherry}
\end{figure}
{By construction, the flow $\phi_-^t$ coincides with the flow $g^t$ for $|x|\geqslant2$. Moreover, the diffeomorphism $\phi_-=\phi_-^1$ has exactly two fixed points: the saddle $P_-(1,0)$ and the sink $Q_-(-1,0)$ (see Figure \ref {cherry}), besides both points are hyperbolic. One unstable separatrix of the $P_-$ saddle is an open interval $(-1,1)\times\{0\}$ belonging to the $Q_-$ sink basin, and the other is a ray $(1,+\infty)\times\{0\}$.

{Define the flow $\phi_+^t$ on $C$ by the formulas
	$$\begin{cases}
		\dot{x}=\begin{cases}1-\frac{1}{9}(x^2+d^2-4)^2, \quad x^2+d^2 \leq 4 \cr
			1, \quad\mbox{else}
		\end{cases}\cr
		\dot{d}=\begin{cases}
			-\frac{d}{2}\big(\sin\big(\frac{\pi}{2}\big(x^2+d^2-3\big)\big)-1\big), \quad 2<x^2+d^2\leq 4\cr
			d,\quad \quad x^2+d^2\leq 2\cr
			0, \quad\mbox{else}
		\end{cases}\cr
	\end{cases}$$} 
\begin{figure}[h]
\centerline{\includegraphics[width=10 true cm, height=7  true cm]{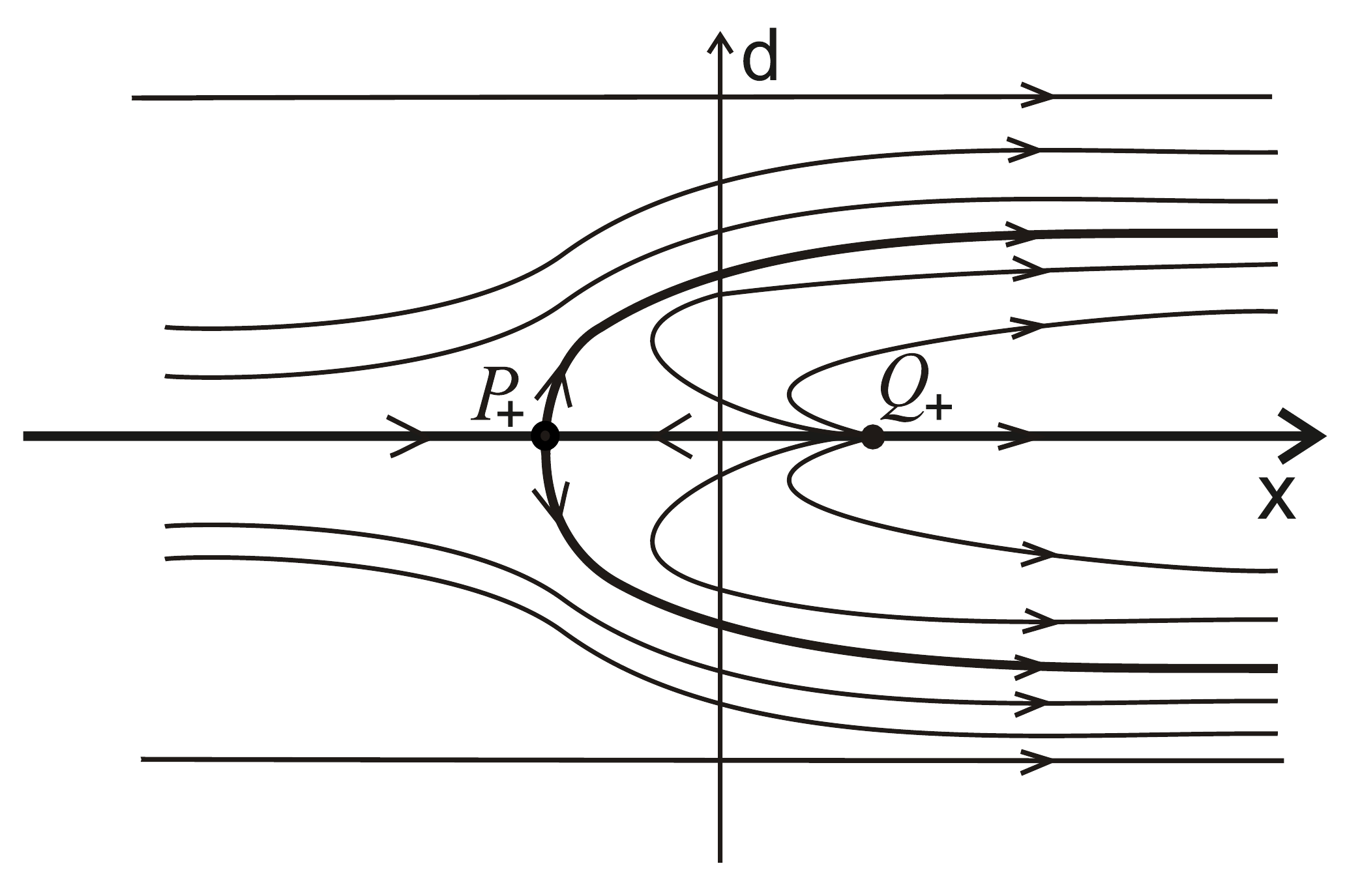}}\caption{Flow paths $\phi_+^t$}\label{cherry2}
\end{figure}
{By construction, the flow $\phi_+^t$ coincides with the flow $g^t$ for $|x|\geqslant 2$. In this case, the diffeomorphism $\phi_+=\phi_+^1$ has exactly two fixed points: the saddle $P_+(-1.0)$ and the source $Q_+(1.0)$ (see Figure \ref{cherry2}), similarly, both points are hyperbolic. One stable separatrix of the saddle $P_+$ is the open interval $(-1,1)\times\{0\}$ belonging to the source basin $Q_+$, and the other is the ray $(-\infty,-1)\times\{0\}$.

Define diffeomorphism $\bar f:\mathbb R^2\to\mathbb R^2$ so that $\bar f$ coincides $h$ outside $A_+\cup A_-$ and coincides $\eta_-^{-1}\phi_-\eta_-$ and $\eta_+^{- 1}\phi_+\eta_+$ to $A_-$ and $A_+$ respectively.

Consider on $\mathbb R^2$ the annulus
$K=\{(x_1,x_2)\in \mathbb R^2: 1\leqslant x_1^2+x_2^2\leqslant 4\}$.Define the function $\nu:[1,2]\to[1,2]$ (see Figure \ref{ris:image1}) by the formula:  $$ \nu (t) = \begin{cases} 1, & t=1, \\ 1+\frac {1} {1+ \exp \left(\frac{\frac{3} {2} -t} {{(t-1)}^2 {\left(t-2 \right)}^2} \right)}, & 1 <t <2, \cr 2, & t= 2. \end{cases} $$
\begin{figure}[h]
\center{\includegraphics[width=0.3\linewidth]{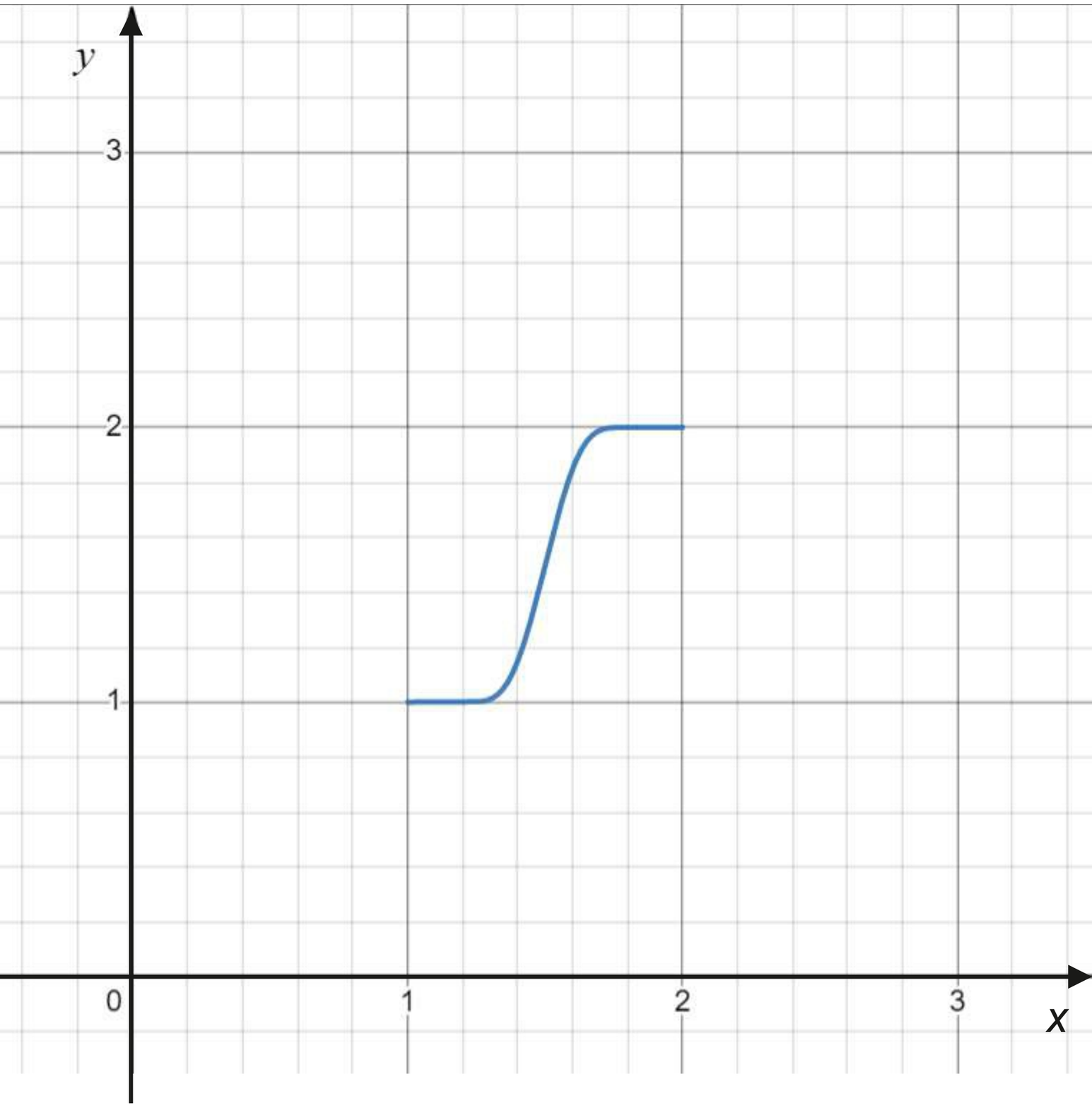}}
	\caption{Function Graph $\nu:[1,2] \to [1,2]$}
	\label{ris:image1}
\end{figure}
On the annulus $K$ we can define the Dehn twist
$\bar d:K \to K$ formula $$\bar d(t,e^{i\phi})=\left(t,e^{i(\phi+ 2\pi \nu(t))}\right).$$
Let $\bar\xi_0=\bar d \circ \bar f:\mathbb R^2 \to \mathbb R^2$ (see Figure \ref{ris:image1}).
\begin{figure}[h]
\center{\includegraphics[width=0.8\linewidth]{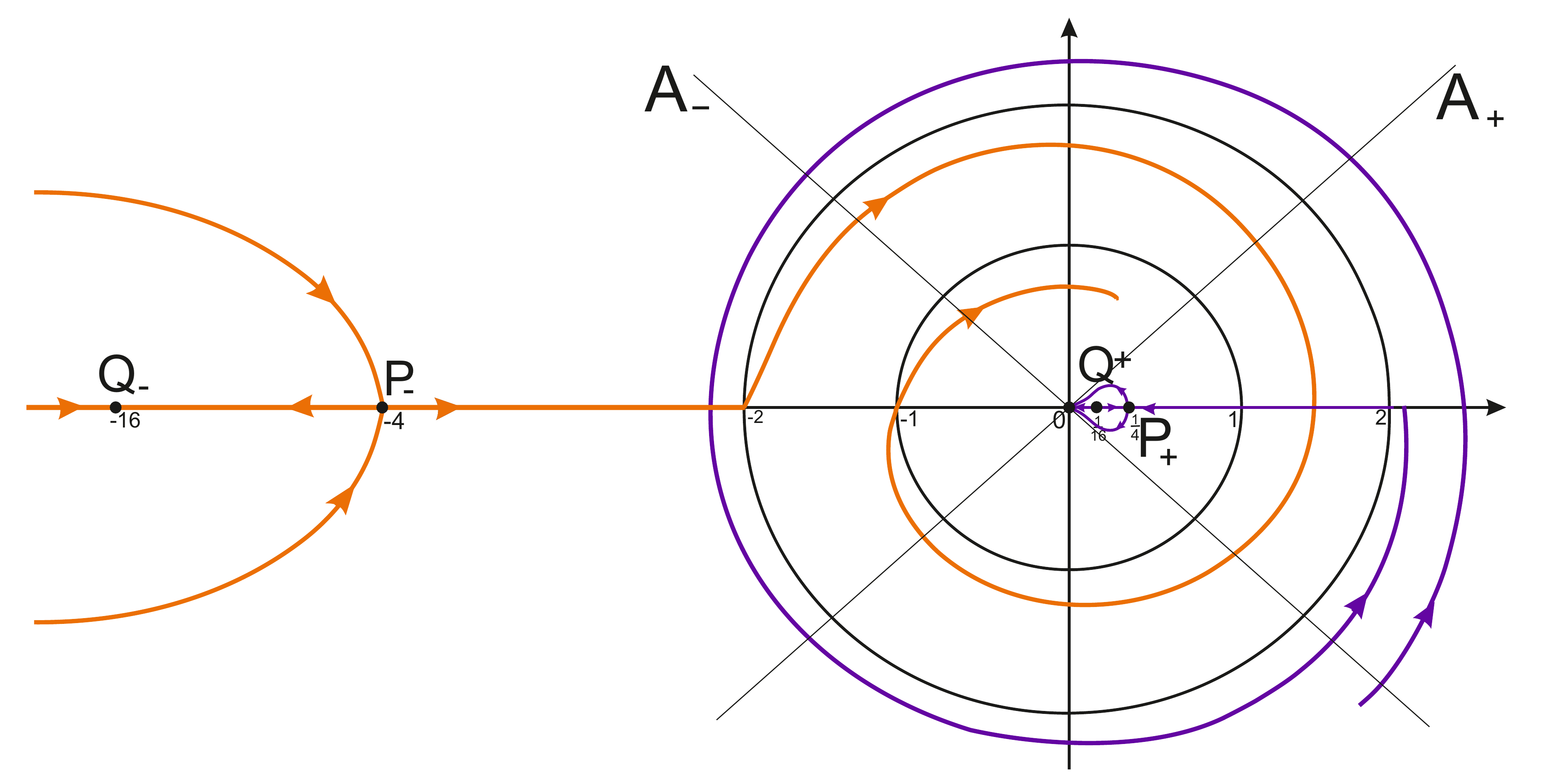}}
	\caption{Diffeomorphism $\bar\xi_0$}
	\label{ris:image1}
\end{figure}

{By construction, the diffeomorphism $\bar\xi_0$ coincides with $h$ in some neighborhood of the point $O$ and the point at infinity, therefore, it induces on $\mathbb{S}^2$ a Morse-Smale diffeomorphism $\xi_0: \mathbb S^2 \to\mathbb S^2$ by the formula $$\xi_0(x)=\begin{cases}
		\vartheta^{-1}\circ \bar \xi_0 \circ \vartheta(x),\,x\in \mathbb S^2\setminus\{N\}, \\
		N,\,x=N.
	\end{cases}$$
It follows directly from the construction that the nonwandering set of the diffeomorphism $\xi_0$ consists of six fixed points of positive orientation: two sources $\alpha_1=N$, $\alpha_2={\xi_0}(\vartheta^{-1}(Q_+ ))$, two sinks $\omega_0={\xi_0}(\vartheta^{-1}(Q_-))$, $\omega=S$ and two saddles $\sigma={\xi_0}(\vartheta^{-1}(P_-)), \sigma_0={\xi_0}(\vartheta^{-1}(P_+))$ (see Figure \ref{net}):
$$\Omega_{\xi_0}=\{\alpha_1,\alpha_2,\omega_0,\omega,\sigma_0,\sigma\}.$$ 
\begin{figure}[h!]
	\center{\includegraphics
		[width=0.6\linewidth]{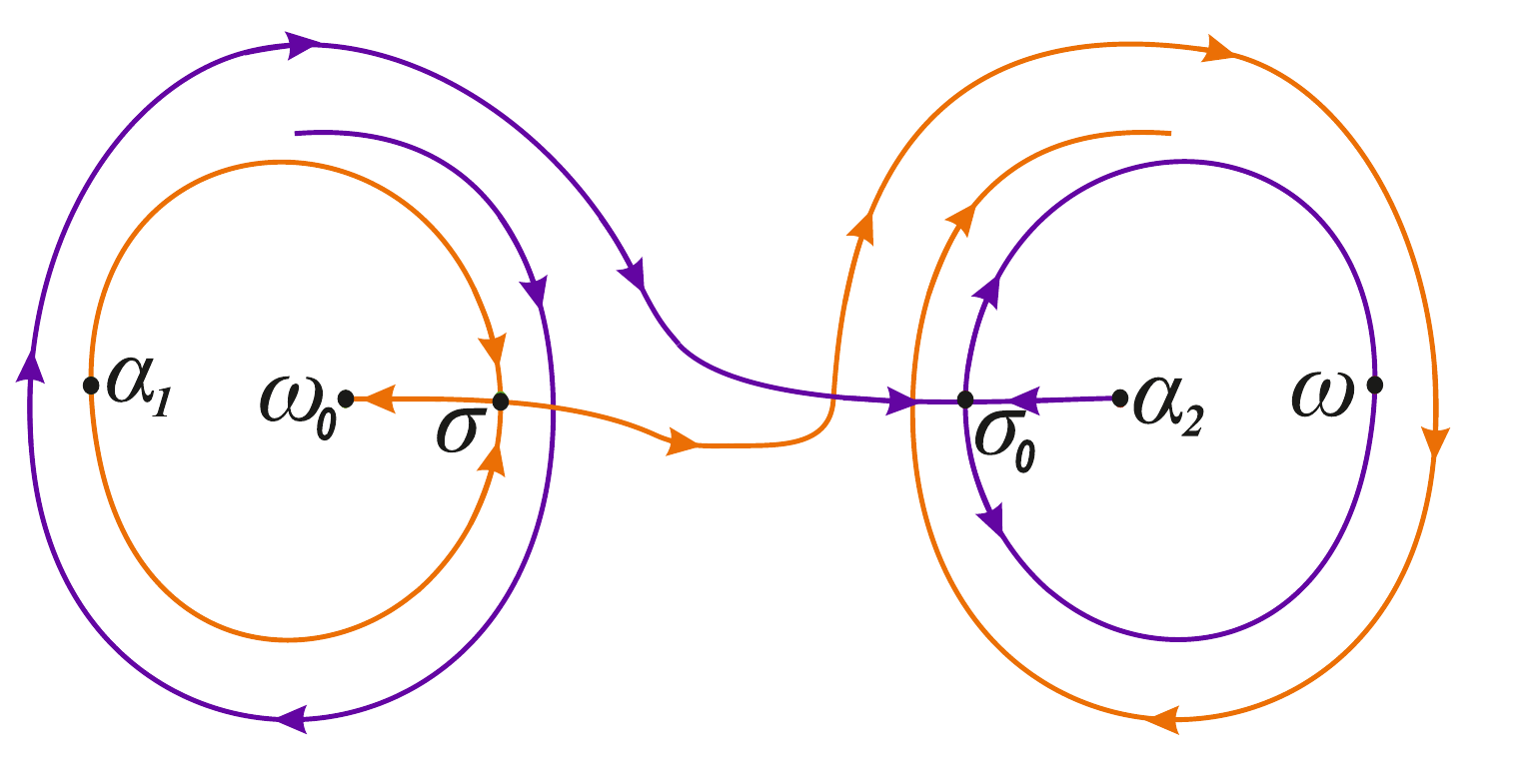}}
	\caption{Phase portrait of a diffeomorphism $\xi_0$}\label{net}
\end{figure}

\section{Proof of the main result}
In this section, we will prove the theorem \ref{th}, each item in a separate lemma below.
\begin{lemma} \label{l41}
On any orientable surface $M^2$ there exists an orientation-changing gradient-like diffeomorphism that does not have a connected characteristic space of orbits.
\end{lemma}
\begin{proof} To prove the lemma, consider a diffeomorphism $f_g:S^+_g\to S^+_g$ such that $f_0=\psi_0$ and $f_g\, (g>0)$ is the connected sum of the diffeomorphism $\psi_0$ with the diffeomorphism $ \psi_g$ along sink $\omega_0$ and source $\alpha$ respectively. By construction, the diffeomorphism $f_g$ changes orientation, and its nonwandering set consists of two fixed sources $\alpha_1,\alpha_2$, two fixed sinks $\omega,\omega_1$ of negative orientation $(\varsigma_{\alpha_1}=\varsigma_ {\alpha_2}=\varsigma_{\omega}=\varsigma_{\omega_1}=-1)$ and $g+1$ saddle periodic orbits $\mathcal O_{\sigma}=\{\sigma,f_g(\sigma )\},\mathcal O_{\sigma_1}=\{\sigma_1,f_g(\sigma_1)\},\dots,\mathcal O_{\sigma_g}=\{\sigma_g,f_g(\sigma_g)\}$  of period 2 and orientation type $(+1,+1)$ (see Figure \ref{SS}):
$$\Omega_{f_g}=\{\alpha, \omega, \sigma,f_g(\sigma_1), \sigma_1,f_g(\sigma_1),\dots, \sigma_{g},f_g(\sigma_g)\}.$$ 
\begin{figure}[h!]
\centerline{\includegraphics
		[width=16 cm]{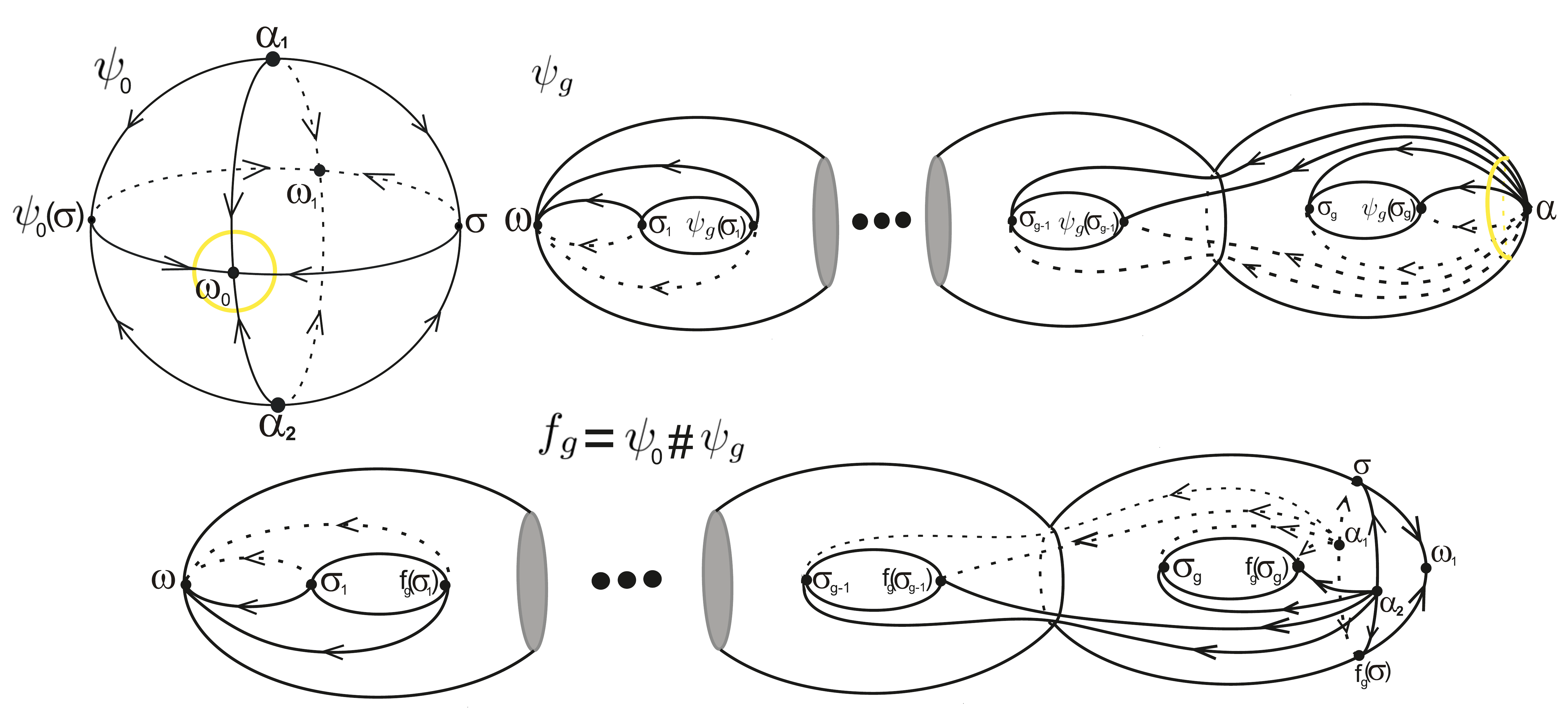}}
	\caption{\small Diffeomorphism $f_g$}\label{Summa_5}
\end{figure}
Let us show that the diffeomorphism $f_g$ does not have a connected characteristic space of orbits.

Indeed, by the proposition \ref{sos1}, each of the orbit spaces $\hat V_{\omega},\,\hat V_{\omega_1}$ is homeomorphic to a Klein bottle. Therefore, if $\Sigma=\emptyset$, then the characteristic orbit space $\hat V_\Sigma$ is not connected and consists of two Klein bottles. Since all saddle points of the diffeomorphism $f_g$ have a positive orientation type, then, according to the corollary \ref{L4}, adding the orbits of such saddles to the set $\Sigma$ does not decrease the number of connected components of the characteristic space of orbits.
\end{proof}

\begin{lemma} On any non-orientable surface $M^2$ there exists a gradient-like diffeomorphism that does not have a connected characteristic space of orbits.
\end{lemma}
\begin{proof}
Define a diffeomorphism $\tilde f_q: S^-_q \rightarrow S^-_q,\,q\in\mathbb N$ as a connected sum of diffeomorphisms $\psi_0$ and $\tilde\psi_q$ ($\tilde f_q=\psi_0 \sharp \tilde\psi_q$) along the sink $\omega_{0}$ of the diffeomorphism $\psi_0$ and the source $\tilde\alpha$ of the diffeomorphism $\tilde\psi_{q}$ (see Fig. \ref{S2} for $q=1$). By construction, the nonwandering set $\Omega_{\tilde f_q}$ consists of two sources $\alpha_1, \alpha_2$ and two sinks $\omega_1, \tilde \omega$ of negative orientation types $(\varsigma_{\omega_1}=\varsigma_ {\tilde\omega}=\varsigma_{\alpha_1}=\varsigma_{\alpha_2}=-1)$, also of $q$ fixed saddles $\tilde\sigma_1,\dots,\tilde\sigma_q$ of orientation type $ \varsigma_{\tilde\sigma_i}=(-1,-1)$ and a saddle orbit $\mathcal O_{\sigma}=\{\sigma,\tilde f_q(\sigma)\}$ of period 2 with a positive orientation type (see Figure \ref{S2} for $q=1$):
$$\Omega_{\tilde f_q}=\{\alpha_1,\alpha_2, \omega_1,\tilde\omega, \sigma,\tilde f_q(\sigma_1),\tilde\sigma_1,\dots,\tilde\sigma_q\}.$$  
\begin{figure}[h!]
\centerline{\includegraphics[width=14 cm]{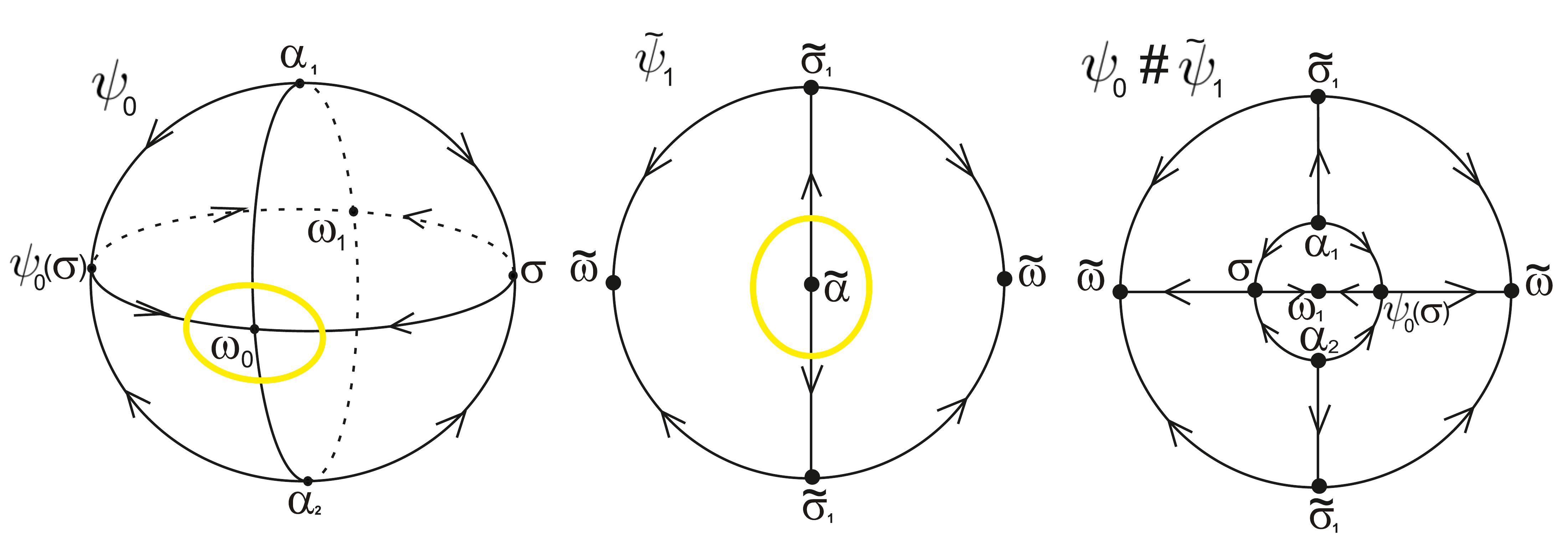}}
	\caption{\small Diffeomorphism $\tilde f_1$}\label{S2}
\end{figure}

Let us show that the diffeomorphism $\tilde f_q$ does not have a connected characteristic space of orbits.

By the construction and proposition \ref{sos1}, each of the orbit spaces $\hat V_{\tilde \omega},\,\hat V_{\omega_1}$ is homeomorphic to a Klein bottle. If $\Sigma=\emptyset$, then the characteristic orbit space $\hat V_\Sigma$ is not connected and consists of two Klein bottles. Since all saddle points of the diffeomorphism $\tilde f_q$ have orientation type either $(+1,+1)$ or $(-1,-1)$, then, according to the corollary \ref{L4}, adding the orbits of such saddles to the set $\Sigma$ does not decrease the number of connected components of the characteristic space of orbits.
\end{proof}

\begin{lemma} \label{l43}
On any surface $M^2$ there exists a Morse-Smale diffeomorphism with heteroclinic points that does not have a connected characteristic space of orbits.
\end{lemma}
\begin{proof}
We construct a diffeomorphism $\xi_g:S^+_g \rightarrow S^+_g$ as a connected sum of diffeomorphisms $\xi_0$ and $\psi^2_g$ $(\xi_g=\xi_0\sharp\psi^2_g)$ along the sink $ \omega_0$ of the diffeomorphism $\xi_0$ and the source $\alpha$ of the diffeomorphism $\psi^2_g$. By construction, the diffeomorphism $\xi_g$ preserves orientation, its nonwandering set consists of points of positive orientation type: two sinks $\omega, \omega_1$ and two sources $\alpha_1, \alpha_2$ and $2+2g$ saddles $\sigma_0,\sigma,\sigma_1,\dots,\sigma_{2g}$ (a special case for $g=1$ is shown in Figure \ref{ex_4}).
\begin{figure}[h!]
\centerline{\includegraphics[width=14 cm]{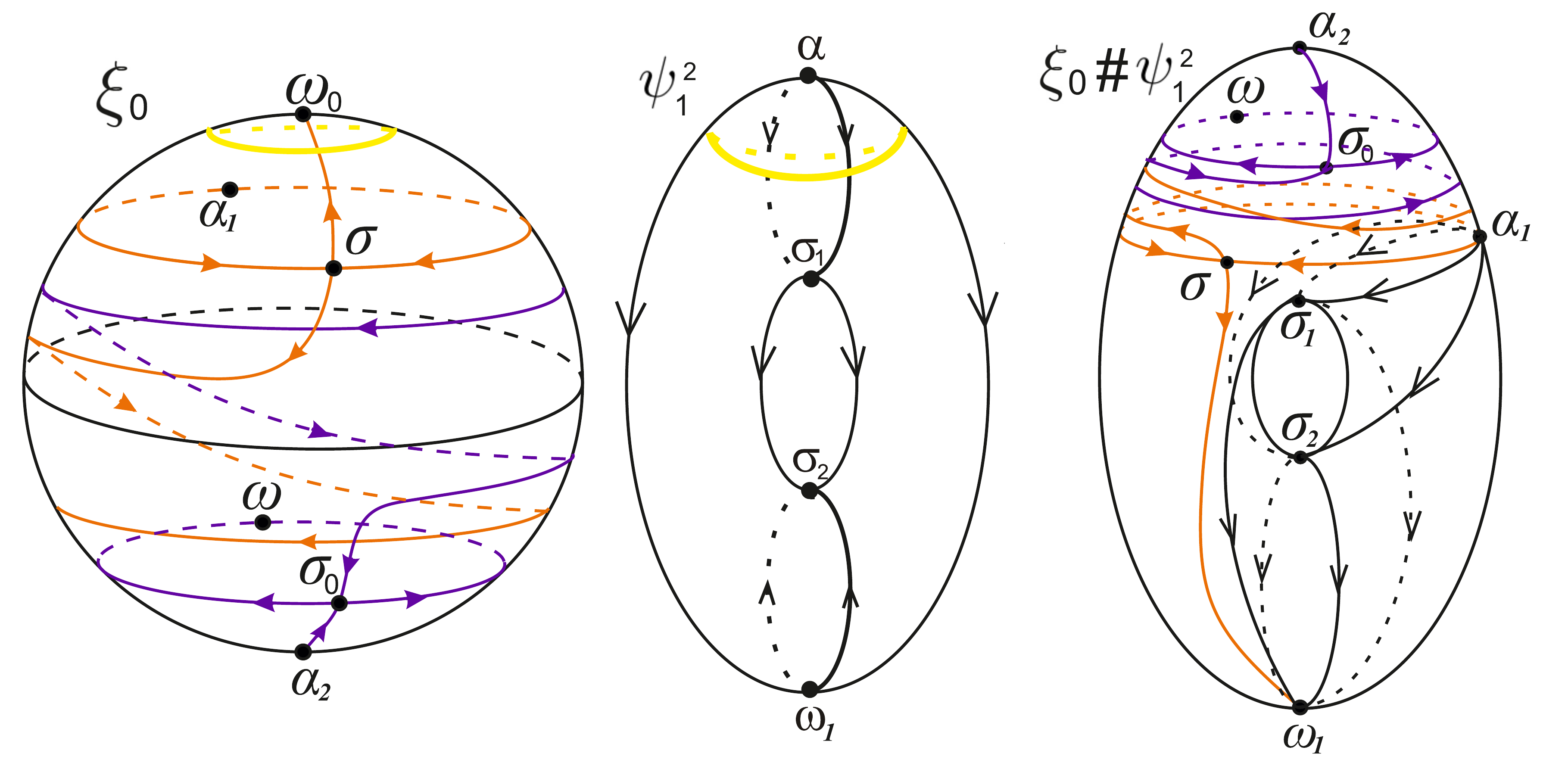}}
	\caption{\small Diffeomorphism $\xi_1$}\label{ex_4}
\end{figure}

Let us show that the diffeomorphism $\xi_g$ does not have a connected characteristic space of orbits.

Let the set $\Sigma=\emptyset$, then the characteristic orbit space consists of two tori $\hat{V}_{\omega_1}$ and $\hat{V}_{\omega}$. Since the unstable saddle point manifolds $\sigma_1,\dots,\sigma_{2g}$ contain the only one sink $\omega$ in their closures, then
adding these saddles to the set $\Sigma$ does not reduce the number of components of the characteristic space of orbits. The unstable manifold of the saddle $\sigma_0$ contains the only one sink $\omega_1$ in its closure, then adding this saddle to the set $\Sigma$ increases the number of connected components of the characteristic space to three. The saddle $\sigma$ is above the saddle $\sigma_0$ by the Smale order, so the saddle $\sigma$ can be added to the set $\Sigma$ only together with the saddle $\sigma_0$. Whence it follows that for any set $\Sigma$ the number of connected components of the characteristic space is greater than one.

ALos, let us construct a diffeomorphism $\tilde\xi_q:S^-_q \rightarrow S^-_q$ as a connected sum of diffeomorphisms $\xi_0$ and $\tilde\psi_q^2$ ($\zeta_q=\xi_0\sharp\tilde\psi_q ^2$) along the sink $\omega_{0}$ of the diffeomorphism $\xi_0$ and the source $\tilde\alpha$ of the diffeomorphism $\tilde\psi_{q}^2$ (a special case for $q=1$ is shown in Figure \ref{ex4}). The nonwandering set of the diffeomorphism $\zeta_q$ consists of points of positive orientation type: two sources $\alpha_1, \alpha_2$, two sinks $\omega_1, \tilde \omega$ and $q+2$ saddles $\sigma, \sigma_0, \tilde \sigma_1, \dots \tilde \sigma_q$.

\begin{figure}[h!]
\centerline{\includegraphics[width=14 cm]{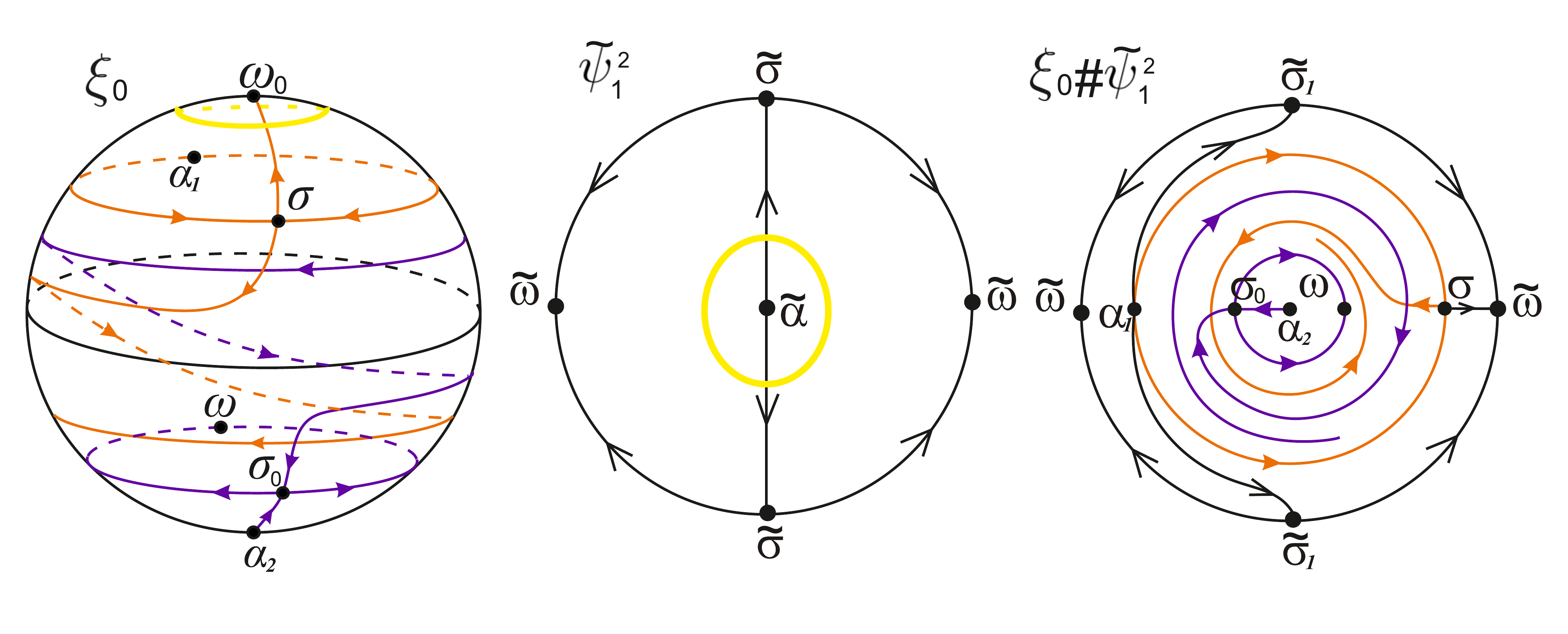}}
	\caption{\small Connected sum of diffeomorphisms $\xi_0$ and $\tilde\psi_1^2$}\label{ex4}
\end{figure}

Similar Arguments to the case of orientable surfaces prove that the diffeomorphism $\tilde\xi_q$ does not have a connected characteristic space of orbits.

\end{proof}

\bibliographystyle{plain} 
\bibliography{refs}

\end{document}